\PassOptionsToPackage{intlimits,tbtags}{amsmath}
\documentclass[aip,jmp,amsmath,amssymb,reprint,onecolumn,longbibliography,showkeys]{revtex4-1}
\usepackage[utf8]{inputenc}
\usepackage[T1]{fontenc}
\usepackage{mathptmx}
\usepackage{amsthm}
\usepackage{microtype}
\usepackage{mathrsfs}
\DeclareFontFamily{U}{rsfs}{\skewchar\font127}
\DeclareFontShape{U}{rsfs}{m}{n}{<-6.5>rsfs5 <6.5-8>rsfs7 <8->rsfs10}{}
\usepackage{enumerate}
\usepackage{etoolbox}
\usepackage[hidelinks]{hyperref}

\hypersetup{
  pdftitle={Rupture and far-field profiles for Henon-type MEMS elliptic equation with pressure},
  pdfauthor={Yunxiao Li and Yanyan Zhang},
  pdfsubject={Manuscript for Journal of Mathematical Physics},
  pdfkeywords={rupture solutions; asymptotic expansions; far-field behavior; singular elliptic equations; weighted nonlinearities; external pressure}
}

\makeatletter
\renewcommand{\p@subsection}{}
\renewcommand{\p@subsubsection}{}
\def\@email#1#2{%
 \endgroup
 \patchcmd{\titleblock@produce}
  {\frontmatter@RRAPformat}
  {\frontmatter@RRAPformat{\produce@RRAP{*#1\href{mailto:#2}{#2}}}\frontmatter@RRAPformat}
  {}{}
}
\makeatother
\AtBeginDocument{\setcitestyle{numbers,square,citesep={], [}}}

\theoremstyle{plain}
\newtheorem{thm}{Theorem}[section]
\newtheorem{cor}[thm]{Corollary}
\newtheorem{lem}[thm]{Lemma}
\newtheorem{prop}[thm]{Proposition}
\theoremstyle{definition}
\newtheorem*{remark}{Remark}
\newtheorem*{defn}{Definition}
\numberwithin{equation}{section}

\begin{document}

\title[Rupture and far-field profiles for H\'enon-type MEMS equations]{Rupture and far-field profiles for H\'enon-type MEMS elliptic equation with pressure}
\author{Yunxiao Li}
\thanks{Electronic mail: \href{mailto:52290155042@stu.ecnu.edu.cn}{52290155042@stu.ecnu.edu.cn}}
\affiliation{School of Mathematical Sciences, East China Normal University,
500 Dongchuan Road, Shanghai 200241, China}
\author{Yanyan Zhang}
\email[Corresponding author: ]{yyzhang@math.ecnu.edu.cn}
\affiliation{School of Mathematical Sciences, Key Laboratory of MEA (Ministry of Education)
and Shanghai Key Laboratory of PMMP, East China Normal University,
500 Dongchuan Road, Shanghai 200241, China}

\begin{abstract}
Motivated by membrane touchdown in electrostatic microelectromechanical
systems (MEMS) under external pressure, we study the H\'enon-type elliptic
equation $\Delta u=\lambda|x|^{\alpha}u^{-p}+F$, $u>0$, where $N\ge2$,
$\lambda,p>0$, $\alpha>-2$, $F\in\mathbb R$, and $\alpha<2p$ if $F\ne0$.
Near the origin, we construct positive radial rupture solutions with
$u(0)=0$ for $F\ne0$ and infinitely many non-radial but asymptotically
radial ones for every $F$, and derive expansions to arbitrary order.
On exterior domains $\Omega_{R_0}:=\{x\in\mathbb R^N:|x|>R_0\}$, with no
condition at the origin, we establish existence and uniqueness of positive
solutions on $[R_0,\infty)$ for radial initial value problems with
$u(R_0)>0$ if $p\ge1$, or if $F\ge0$ and $u'(R_0)\ge0$.
We determine the leading behavior of all positive radial exterior solutions
as $|x|\to\infty$ for $F>0$ and $F=0$, and for $F<0$ when $p\ge1$ and
$\alpha\ge0$.
For possibly non-radial exterior solutions satisfying
$u(r,\theta)=A_\infty r^{\sigma_\infty}(1+o(1))$ uniformly in $\theta$,
with $A_\infty>0$ and $\sigma_\infty\in\mathbb R$, we uniquely determine the leading exponent and coefficient
in all three pressure regimes, without derivative bounds on the remainder.
Although pressure affects only higher-order terms near an asymptotically
radial rupture point, its sign changes the leading-order behavior at
infinity under the stated hypotheses.
Local radial rupture solutions extend uniquely to
$\mathbb R^N\setminus\{0\}$ for $F>0$, and for $F<0$ when $p\ge1$,
with the classified far-field profiles under the corresponding hypotheses.
\end{abstract}


\maketitle

\section{Introduction}

We study positive solutions of the elliptic equation
\begin{equation}\label{QQ}
 \Delta u=\frac{\lambda |x|^{\alpha}}{u^p}+F,
 \qquad u>0,
\end{equation}
where $N\ge2$ is the space dimension, $\lambda>0$, $p>0$, $\alpha>-2$, and
$F\in\mathbb R$.  Equation \eqref{QQ} generalizes the stationary equation
of an electrostatic membrane model for microelectromechanical systems
(MEMS) with external pressure; see \cite{23z}, \cite{24z}, \cite{25z}
and \cite{26z} for related models.  In the classical interpretation,
$u$ is the normalized distance between a deformable membrane and a fixed
ground plate, $\lambda$ is a parameter determined by the applied voltage,
the H\'enon weight $|x|^\alpha$ describes an idealized spatially varying
permittivity of the membrane, and $F$ is an external pressure.
Contact between the membrane and the plate corresponds to the vanishing
of $u$ and is called touchdown, or rupture; accordingly, a positive
solution of \eqref{QQ} in a punctured ball $B_R\setminus\{0\}$ with
$u(0)=0$ is called a rupture solution with an isolated rupture point at
the origin.  In this stationary setting, the vanishing rate and
higher-order terms of $u$ describe the local membrane shape near contact.
The classical planar model has $p=2$ and $N=2$; we allow general $p>0$,
any dimension $N\ge2$ and pressures of either sign, and we make explicit
how the profiles at a rupture point and at infinity depend on
$\alpha$, $F$, $p$ and $N$.

Throughout the paper we assume
\begin{equation*}\tag{H}\label{H}
 N\ge2,\qquad \lambda>0,\qquad p>0,\qquad \alpha>-2,\qquad
 F\in\mathbb R,\qquad \alpha<2p\quad\text{if }F\ne0,
\end{equation*}
unless the contrary is stated explicitly.

A basic role is played by the power law obtained by balancing the Laplacian
against the singular H\'enon term.  Set
\begin{equation}\label{qcc}
 d:=\frac{\alpha+2}{p+1},
 \qquad
 \Lambda:=\left(\frac{\lambda}{d(d+N-2)}\right)^{\frac1{p+1}}.
\end{equation}
Since $N\ge2$ and $\alpha>-2$, both $d$ and $\Lambda$ are positive.  For
$r=|x|$ and $C>0$ one has
\[
 \Delta(Cr^d)=Cd(d+N-2)r^{d-2},
 \qquad
 \frac{\lambda r^\alpha}{(Cr^d)^p}=\lambda C^{-p}r^{\alpha-pd},
\]
and matching exponents and coefficients gives $d=(\alpha+2)/(p+1)$ and
$C=\Lambda$.  Thus $u(x)=\Lambda|x|^d$ is an exact solution of \eqref{QQ}
when $F=0$; we call it the balanced profile.  When $F\ne0$, writing
$u(r,\theta)=r^{d}\bigl(\Lambda+z(\log r,\theta)\bigr)$ and dividing \eqref{QQ} by $r^{d-2}$
turns the constant $F$ into $Fr^{2-d}$.  We therefore set
\begin{equation}\label{eq:pressure-exponent-intro}
 \rho:=2-d=\frac{2p-\alpha}{p+1},
\end{equation}
which is positive precisely because $\alpha<2p$.  

\begin{defn}
A positive solution $u$ of \eqref{QQ} in a punctured ball is called a
\emph{positive asymptotically radial rupture solution} if
\begin{equation}\label{defn}
 u(x)=\Lambda|x|^{d}\bigl(1+O(|x|^\varepsilon)\bigr)
 \qquad\text{as } |x|\to0^{+}
\end{equation}
for some $\varepsilon>0$, where $\Lambda$ and $d$ are given by \eqref{qcc}.
Here $f(x)=O(g(x))$ as $x\to x_0$ means that $|f(x)|\le C|g(x)|$ for some
constant $C>0$ and all $x$ sufficiently close to $x_0$.
\end{defn}

Two asymptotic problems are addressed.  The first is local: in a punctured
ball we ask whether rupture solutions of \eqref{QQ} exist for every
pressure, and we determine the complete asymptotic behavior of positive
asymptotically radial rupture solutions at the rupture point.  The second is
global: in exterior domains, with no condition imposed at the origin,
we establish radial existence by continuation from positive initial data
and classify the possible leading-order behaviors of positive solutions
as $|x|\to\infty$. This exterior analysis concerns the equation on
unbounded domains, rather than a boundary-value problem for a finite device.
The two problems are independent of each other: a solution of the first need not extend beyond a small ball, and a solution of the second need not extend to the origin. In the radial setting, however, the two problems can be linked for
nonzero pressure.  The local rupture solutions of
Theorem~\ref{thm:radial} extend uniquely to
$\mathbb R^N\setminus\{0\}$ when $F>0$, and also when $F<0$ with
$p\ge1$; in the parameter ranges covered by the far-field classification,
they connect the rupture profile at the origin to the corresponding
profile at infinity (Corollary~\ref{cor:radial-global-connection}).

We briefly recall the results most closely related to ours.

\smallskip\noindent
\emph{Rupture solutions for $\alpha=0$.}  For the unweighted equation,
rupture solutions were studied in \cite{19z} and \cite{21z}.  In particular, for $\lambda=1$, $p>1$ and $N\geq2$, \cite{21z} studies the dimension of the rupture set $\{u=0\}$ in a bounded domain for $F=0$. For $\lambda=1$, $F=0$ and $3\leq N \leq10$, \cite{19z} obtained infinitely many non-radial rupture solutions for $p > \frac{N-1-2\sqrt{N-2}}{2\sqrt{N-2}-(N-5)}$.

\smallskip\noindent
\emph{The H\'enon case $\alpha\ne0$.}  For $\lambda=1$, $F=0$, $N\ge2$,
$p>0$ and $\alpha>-2$, \cite[Theorem~2.3]{12z} proved that the balanced
profile is the unique positive radial solution of \eqref{QQ} in a punctured
ball with an isolated rupture at the origin, that it extends to all $r>0$; the statements for
general $\lambda>0$ follow from the elementary amplitude rescaling that
normalizes $\lambda$ to one.  The two-dimensional problem was studied in
\cite{7z} and \cite{9z}: for $p=2$ and $F\ge0$, \cite{7z} identified parameter ranges
in which the leading rupture profile is isotropic and ranges in which it may
be anisotropic; for $F=0$ and general $p>0$, \cite{9z} described the
behavior of rupture solutions, including global solutions on
$\mathbb R^2\setminus\{0\}$. 

\smallskip\noindent
\emph{Expansions near isolated singularities.} Some research has also been made regarding the asymptotic expansion of solutions near singular points. For instance, \cite{4z} derived the asymptotic expansion near an isolated singularity for the Yamabe equation $-\Delta u = \frac{1}{4}n(n-2)u^{\frac{n+2}{n-2}}$ by using spherical harmonics. Then \cite{8z} analyze asymptotic expansions to arbitrary order for the steady thin-film-type equation near the origin inspired by \cite{4z}, which corresponds to a special case of Equation \eqref{QQ} with $\lambda=1$ and $\alpha=0$. For $\alpha>-2$, \cite{7z} derived the first two terms of the expansion for asymptotically isotropic rupture solutions in the case of $N=2$. 

\smallskip\noindent
\emph{Behavior at infinity.}  Much less is known about the far field. For
$F=0$, Du and Guo~\cite[Theorem~2.3 and Lemma~3.2(i)]{12z} showed that the
balanced profile is a global radial rupture solution and constructed a
family of positive radial entire solutions asymptotic to this profile at
infinity.  Global rupture solutions on $\mathbb R^2\setminus\{0\}$ are
studied in \cite{9z}, and the behavior at infinity of entire radial solutions
of $\Delta u=\mathcal W(|x|)u^{-p}$ with bounded nonincreasing $\mathcal W$ is determined in
\cite{6z}.  For $F<0$ and $\alpha=0$,
Jiang and Ni \cite{27z} ($N\ge2$, $p>1$) and Guo, Ye and Zhou \cite{28z} showed,
 under their respective assumptions, that nonconstant entire radial solutions oscillate about and converge to the equilibrium  $(\lambda/(-F))^{1/p}$. 

We now state the main results.  The first concerns radial rupture solutions
for nonzero pressure.

\begin{thm}\label{thm:radial}
Assume \eqref{H} and $F\ne0$.  Then there exist $R>0$ and a positive
radial rupture solution $u\in C^2(B_R\setminus\{0\})$ of \eqref{QQ}
such that
\begin{equation}\label{eq:radial-origin-asymptotic}
 u(r)=\Lambda r^d\bigl(1+O(r^\rho)\bigr)
 \qquad\text{as }r\to0^+,
\end{equation}
where $r=|x|$, $\Lambda$ and $d$ are given by \eqref{qcc}, and $\rho$ by
\eqref{eq:pressure-exponent-intro}.  Moreover, for every integer $K\ge1$,
\[
 u(r)=\Lambda r^d+\sum_{i=1}^K d_i r^{d+i\rho}
 +O\bigl(r^{d+(K+1)\rho}\bigr)
 \qquad\text{as }r\to0^+.
\]
The coefficients satisfy
$d_i=F^iD_i(\lambda,p,\alpha,N)$ and are determined recursively in
Proposition~\ref{prop:radial-exp}.
\end{thm}

\begin{remark}
Theorem~\ref{thm:radial} is stated for $F\ne0$ because for $F=0$ the
balanced profile $\Lambda r^d$ is the unique positive radial rupture
solution in a punctured ball, as recalled above. 

\end{remark}

To state the non-radial results we introduce the spectral notation.  Let
\begin{equation}\label{eq:AB}
 A_0=N-2+2d>0,
 \qquad
 B_0=(\alpha+2)(N-2+d)>0,
 \qquad
 \lambda_k=k(N-2+k),
\end{equation}
where $\lambda_k$ is the eigenvalue of $-\Delta_{S^{N-1}}$ associated
with spherical harmonics of degree $k$.  Define the modal indicial
polynomials
\begin{equation}\label{eq:modal-indicial-polynomial}
 P_k(s):=s^2+A_0s+B_0-\lambda_k,
 \qquad k=0,1,2,\ldots,
\end{equation}
and set $\Delta_k=A_0^2-4(B_0-\lambda_k)$.  Their roots are
\begin{equation}\label{eq:modal-roots}
 \sigma_1^{(k)}=\frac{-A_0+\sqrt{\Delta_k}}2,
 \qquad
 \sigma_2^{(k)}=\frac{-A_0-\sqrt{\Delta_k}}2,
\end{equation}
with the usual complex-conjugate interpretation when $\Delta_k<0$.
Since $\lambda_k$ is strictly increasing and tends to infinity,
\begin{equation}\label{eq:kstar}
 k_*:=\min\{k\ge1:\lambda_k>B_0\}
\end{equation}
is well defined.  Equivalently, $k_*$ is the unique positive integer such
that $\delta^{(k_*-1)}\le\alpha<\delta^{(k_*)}$, where
\begin{equation}\label{eq:delta}
\delta^{(k)}=
\begin{cases}
-2, & k=0,\\[3pt]
-\dfrac12(p+1)(N-2)
+\dfrac12\sqrt{(N-2)^2(p+1)^2+4(p+1)k(N-2+k)}-2,
& k\ge1.
\end{cases}
\end{equation}
For $k\ge k_*$ one has $\sigma_1^{(k)}>0>\sigma_2^{(k)}$; write
$\gamma_k:=\sigma_1^{(k)}$.  The sequence $\{\gamma_k\}_{k\ge k_*}$ is
strictly increasing and tends to infinity, whereas for $k<k_*$ both
indicial roots have nonpositive real part.  The generating exponents are
\begin{equation}\label{eq:GF}
 \mathcal G_F:=\{\gamma_k:k\ge k_*\}\cup
 \begin{cases}
  \{\rho\}, & F\ne0,\\[2pt]
  \emptyset, & F=0,
 \end{cases}
\end{equation}
and $\mathcal I$ denotes the positive additive semigroup generated by
$\mathcal G_F$:
\begin{equation}\label{eq:I-sets}
 \mathcal I:=
 \left\{
  \sum_{\gamma\in\mathcal G_F}n_\gamma\gamma:\quad
  n_\gamma\in\mathbb N_0\ \text{finitely supported},\quad
  \sum_{\gamma\in\mathcal G_F}n_\gamma\ge1
 \right\}.
\end{equation}
Because $\mathcal G_F$ has a positive minimum and no finite accumulation
point, $\mathcal I$ is locally finite and can be written as
\begin{equation}\label{eq:mu-seq}
 \mu_1<\mu_2<\mu_3<\cdots\to\infty,
 \qquad
 \mu_1=\min\mathcal G_F=
 \begin{cases}
  \sigma_1^{(k_*)}, & F=0,\\[2pt]
  \min\{\sigma_1^{(k_*)},\rho\}, & F\ne0.
 \end{cases}
\end{equation}
These exponents depend only on $N,p,\alpha$ and on whether $F$ vanishes;
they are independent of the particular solution.

\begin{thm}\label{thm:nonradial}
Assume \eqref{H}, and let $d$, $\Lambda$, $\rho$, $k_*$,
$\sigma_1^{(k_*)}$ and $\{\mu_j\}_{j\ge1}$ be as in
\eqref{qcc}--\eqref{eq:pressure-exponent-intro} and
\eqref{eq:AB}--\eqref{eq:mu-seq}.
\begin{enumerate}[(i)]
\item \emph{(Existence of non-radial asymptotically radial solutions.)}
There exist $R>0$ and infinitely many positive asymptotically radial
rupture solutions $u\in C^2(B_R\setminus\{0\})$ of \eqref{QQ} which are
not radially symmetric.
\item \emph{(Asymptotic expansion.)} Let $u\in C^2(B_{R'}\setminus\{0\})$,
for some $R'>0$, be any positive asymptotically radial rupture solution of
\eqref{QQ}.  Then
$$
u(x)=\Lambda |x|^{d}\left(1+O\left(|x|^{\mu_1}\right)\right)
\qquad \text{as } |x| \rightarrow 0^{+},
$$
where $d$ is given by \eqref{qcc}, $\mu_1=\sigma_1^{(k_*)}$ if $F=0$ and
$\mu_1=\min\{\sigma_1^{(k_*)},\rho\}$ if $F\ne0$.  Moreover, for every
integer $K\ge1$,
$$
u(x)=\Lambda |x|^{d}+\sum_{j=1}^{K}\sum_{i=0}^{M_j}
 C_{ji}\!\left(\frac{x}{|x|}\right)(\ln |x|)^i|x|^{d+\mu_j}
+O\!\left((1+|\ln|x||)^{M_{K+1}}|x|^{d+\mu_{K+1}}\right)
$$
as $|x|\to0^+$, where $M_1,\ldots,M_{K+1}$ are nonnegative integers and
each $C_{ji}$ is a finite sum of spherical harmonics on $S^{N-1}$.  The
exponents $\mu_j$ are the elements of the semigroup $\mathcal I$ in
\eqref{eq:I-sets}--\eqref{eq:mu-seq} and depend only on $N$, $p$, $\alpha$
and on whether $F=0$, whereas the integers $M_j$ and the coefficients
$C_{ji}$ may depend on $u$.
\end{enumerate}
\end{thm}

We next turn to the far-field problem.  For $R_0>0$ let
$\Omega_{R_0}:=\{x\in\mathbb R^N:|x|>R_0\}$, and consider positive classical
solutions $u\in C^2(\Omega_{R_0})$ of
\begin{equation}\label{eq:infinity-equation}
 \Delta u=\frac{\lambda|x|^\alpha}{u^p}+F
 \qquad\text{in }\Omega_{R_0},
\end{equation}
under the standing assumption \eqref{H}; no condition is imposed at the
origin.  For possibly non-radial solutions we assume only a uniform
power-type asymptotic: writing $r=|x|$ and $\theta=x/|x|$,
\begin{equation}\label{eq:infinity-power-hypothesis}
 u(r,\theta)=A_\infty r^{\sigma_\infty}\bigl(1+o(1)\bigr),
 \qquad A_\infty>0,\quad \sigma_\infty\in\mathbb R,
 \qquad r\to\infty,
\end{equation}
with the remainder uniform in $\theta\in S^{N-1}$ and no assumption on its
derivatives.

\begin{thm}[Far-field classification]\label{thm:infinity-classification}
Assume \eqref{H}.

\noindent\emph{(1) Radial case.}
For every positive radial solution $u$ of \eqref{eq:infinity-equation}
on $\Omega_{R_0}$, the following asymptotic formulas hold as $r\to\infty$
in the indicated parameter regimes:
\begin{equation}\label{eq:infinity-radial-conclusion}
 u(r)=
 \begin{cases}
  \displaystyle \frac{F}{2N}r^2\bigl(1+o(1)\bigr),
    & F>0,\\[5pt]
  \displaystyle \Lambda r^d\bigl(1+o(1)\bigr),
    & F=0,\\[5pt]
  \displaystyle \left(\frac{\lambda}{-F}\right)^{1/p}
  r^{\alpha/p}\bigl(1+o(1)\bigr),
    & F<0,\ p\ge1,\ \alpha\ge0.
 \end{cases}
\end{equation}
where $d$ and $\Lambda$ are given by \eqref{qcc}.

\noindent\emph{(2) Possibly non-radial case.}
If $u$ is a positive solution of \eqref{eq:infinity-equation}
satisfying \eqref{eq:infinity-power-hypothesis}, then, as $|x|\to\infty$,
\begin{equation}\label{eq:infinity-nonradial-conclusion}
 u(x)=
 \begin{cases}
  \displaystyle \frac{F}{2N}|x|^2\bigl(1+o(1)\bigr),
    & F>0,\\[5pt]
  \displaystyle \Lambda|x|^d\bigl(1+o(1)\bigr),
    & F=0,\\[5pt]
  \displaystyle \left(\frac{\lambda}{-F}\right)^{1/p}
  |x|^{\alpha/p}\bigl(1+o(1)\bigr),
    & F<0.
 \end{cases}
\end{equation}
The remainders are uniform in $\theta$.
\end{thm}

\begin{remark}
Assume \eqref{H}.  For $R_0>0$, $u_0>0$ and $u_1\in\mathbb R$, consider
the radial initial value problem
\begin{equation}\label{eq:radial-initial-value-problem}
 \left\{
 \begin{aligned}
  u''+\frac{N-1}{r}u' &= \lambda r^\alpha u^{-p}+F,
  && r>R_0,\\[3pt]
  u(R_0) &= u_0,\\[3pt]
  u'(R_0) &= u_1.
 \end{aligned}
 \right.
\end{equation}
As proved in Lemma~\ref{lem:radial-continuation}, this problem has a
unique positive solution on $[R_0,\infty)$ if $p\ge1$, or if $F\ge0$
and $u_1\ge0$.
\end{remark}

Applying this result of remark to the local radial rupture solutions
of Theorem~\ref{thm:radial} and combining it with the far-field
classification \eqref{eq:infinity-radial-conclusion} yields the following
corollary.

\begin{cor}[Global radial connections for nonzero pressure]
\label{cor:radial-global-connection}
Assume \eqref{H} and $F\ne0$.  Let
$u\in C^2(B_R\setminus\{0\})$ be a positive radial solution of \eqref{QQ}
satisfying \eqref{eq:radial-origin-asymptotic}.
\begin{enumerate}[(i)]
\item If $F>0$, then $u$ extends uniquely to a positive radial solution on
$\mathbb R^N\setminus\{0\}$ and
\[
 u(r)=\frac{F}{2N}r^2\bigl(1+o(1)\bigr)
 \qquad\text{as }r\to\infty.
\]

\item If $F<0$ and $p\ge1$, then $u$ extends uniquely to a positive radial
solution on $\mathbb R^N\setminus\{0\}$.  If, in addition, $\alpha\ge0$,
then
\[
 u(r)=\left(\frac{\lambda}{-F}\right)^{1/p}
 r^{\alpha/p}\bigl(1+o(1)\bigr)
 \qquad\text{as }r\to\infty.
\]
\end{enumerate}
In particular, Theorem~\ref{thm:radial} supplies global radial rupture
solutions in both cases.
\end{cor}
\begin{remark}

The continuation result in Lemma~\ref{lem:radial-continuation}
imposes no condition at the origin and does not by itself produce
rupture solutions.  Corollary~\ref{cor:radial-global-connection} applies
it to the local radial rupture solutions of Theorem~\ref{thm:radial}.
For $F=0$, the global radial rupture solution is the explicit balanced
profile \cite[Theorem~2.3]{12z}.

\end{remark}

\smallskip\noindent
\textbf{Novelty and main contributions.}
The main contributions of these results, in comparison with the works
reviewed above, are as follows.

\smallskip\noindent
\textbf{(A) Arbitrary-order rupture expansions for H\'enon weights and
external pressure.}
Theorem~\ref{thm:nonradial}(ii) gives an expansion to any prescribed
order, with a quantified remainder, for every positive asymptotically
radial rupture solution.  This extends the arbitrary-order unweighted
analysis in \cite{8z} and the two-term planar analysis in \cite{7z}
to all dimensions $N\ge2$ and the full parameter range \eqref{H}.
The construction also covers the exceptional cases of zero indicial
roots and resonances requiring logarithmic correction terms.
The local existence results in Theorems~\ref{thm:radial}
and~\ref{thm:nonradial}(i) ensure that this expansion theory applies
both to radial rupture solutions for $F\ne0$ and to infinitely many
genuinely non-radial but asymptotically radial rupture solutions for every
$F\in\mathbb R$. In the classical MEMS setting, these expansions describe
the membrane shape near an isolated contact point, including
pressure-dependent and angular corrections to the leading radial profile.

\smallskip\noindent
\textbf{(B) Radial exterior existence and a three-regime far-field classification.}
For positive radial exterior solutions,
Theorem~\ref{thm:infinity-classification}(1) classifies the leading
behavior according to the sign of $F$, without any condition at the
origin; for $F<0$, this conclusion requires $p\ge1$ and
$0\le\alpha<2p$.  Lemma~\ref{lem:radial-continuation} also gives
existence and uniqueness for \eqref{eq:radial-initial-value-problem}
on $[R_0,\infty)$ when $p\ge1$, or when $F\ge0$ and $u_1\ge0$.
For possibly non-radial solutions,
Theorem~\ref{thm:infinity-classification}(2) determines the exponent
and leading coefficient under the uniform power-law hypothesis
\eqref{eq:infinity-power-hypothesis}, without derivative assumptions
on the remainder.  The earlier far-field results reviewed above concern
$F=0$ or the unweighted case $F<0$, $\alpha=0$.  Under the respective
hypotheses, our classification also covers $F>0$ and negative pressure
with nontrivial H\'enon weights. In contrast to its higher-order role
near an asymptotically radial rupture point, pressure changes the
leading-order far-field balance. Positive pressure produces quadratic
growth, zero pressure gives the balanced profile, and negative pressure
is balanced by the weighted nonlinear term $\lambda |x|^\alpha u^{-p}$ at leading order. For nonzero pressure,
its magnitude also enters the leading coefficient.

\smallskip\noindent
\textbf{(C) A local-to-global connection between the two asymptotic
problems.}
Corollary~\ref{cor:radial-global-connection} gives unique global
continuations of the local radial rupture solutions for $F>0$, and
for $F<0$ with $p\ge1$.  When the radial far-field classification
also applies, these solutions connect the rupture expansion at the
origin to the corresponding leading profile at infinity.

Sections~2 and~3 derive arbitrary-order asymptotic expansions near the
rupture point for radial and asymptotically radial solutions, respectively.
Section~4 constructs local radial and non-radial rupture solutions and
completes the proofs of Theorems~\ref{thm:radial} and~\ref{thm:nonradial}.
Section~5 proves, in order, a radial continuation lemma, the far-field
classification (Theorem~\ref{thm:infinity-classification}), and the
rupture-to-infinity connection (Corollary~\ref{cor:radial-global-connection}).

\section{Asymptotic behavior of radial solutions}

Throughout Sections~2--4, set $r=|x|$, $t=\log r$,
$\theta=x/|x|$, and
\[
 z(t,\theta):=r^{-d}u(x)-\Lambda.
\]
In polar coordinates, equation \eqref{QQ} is
\begin{equation}\label{Qz}
 u_{rr}+\frac{N-1}{r}u_r+\frac1{r^2}\Delta_{S^{N-1}}u
 =\frac{\lambda r^\alpha}{u^p}+F.
\end{equation}
Using \eqref{qcc}, a direct computation gives
\begin{equation}\label{3}
 z_{tt}+(N-2+2d)z_t+(\alpha+2)(N-2+d)z
 +\Delta_{S^{N-1}}z
 =f(z)+Fe^{(2-d)t},
\end{equation}
where
\[
 f(z):=\frac{\lambda}{(z+\Lambda)^p}
 -\frac{\lambda}{\Lambda^p}
 +\frac{\lambda p}{\Lambda^{p+1}}z.
\]
For radial solutions, $z$ is independent of $\theta$, and we derive the radial expansion to every finite order.
\begin{prop}\label{prop:radial-exp}
Assume \eqref{H} and $F\ne0$, and put $\rho=\frac{2p-\alpha}{p+1}>0$. Let $u$ be a positive radial rupture solution of \eqref{QQ} such that, for some $\varepsilon>0$,
\[
u(r)=\Lambda r^{\frac{\alpha+2}{p+1}}\left(1+O\left(r^{\varepsilon}\right)\right)\qquad\text{as } r \rightarrow 0^{+},
\]
where $r=|x|$ and $\Lambda$ is given by \eqref{qcc}, and set
$z(t):=r^{-d}u(r)-\Lambda$ with $t=\ln r$. Then, for every integer $K\ge1$,
\begin{equation*}
z(t) = \sum_{\ell=1}^{K} d_{\ell} e^{\ell \rho t} + O \left( e^{(K+1)\rho t} \right)\qquad\text{as } t\to-\infty,
\end{equation*}
where, with $P_0(s)=s^2+A_0s+B_0$ the polynomial \eqref{eq:modal-indicial-polynomial} for $k=0$ (so that $A_0$, $B_0$ are given by \eqref{eq:AB}) and
\[
a_m=\frac{\lambda(-1)^m(p)_m}{m!\,\Lambda^{p+m}},\qquad
(p)_m=p(p+1)\cdots(p+m-1),
\]
the coefficients are determined recursively by
\begin{equation}\label{eq:radial-recurrence}
d_1=\frac{F}{P_0(\rho)},\qquad
d_n=\frac{1}{P_0(n\rho)}\sum_{m=2}^{n}a_m
\sum_{\substack{j_1+\cdots+j_m=n\\j_1,\ldots,j_m\ge1}}
d_{j_1}\cdots d_{j_m}\quad(n\ge2).
\end{equation}
In particular, $d_n=F^nD_n(\lambda,p,\alpha,N)$.
\end{prop}

\begin{proof}
For radial $u$, equation \eqref{3} reduces to
\begin{equation}\label{eq:radial-ode}
 P_0(\partial_t)z=f(z)+Fe^{\rho t},\qquad
 z(t)=O(e^{\varepsilon t})\quad\text{as }t\to-\infty,
\end{equation}
where
\[
 P_0(s)=s^2+A_0s+B_0,\qquad
 A_0=N-2+2d,\qquad B_0=(p+1)d(N-2+d).
\]
These coefficients agree with \eqref{eq:AB}, since $(p+1)d=\alpha+2$.
The function $f$ is smooth for $|s|\le\Lambda/2$ and satisfies
\[
 f(0)=f'(0)=0,\qquad f(s)=O(s^2),\qquad f'(s)=O(s).
\]
Fix $T<\min\{0,\ln R\}$ sufficiently negative that
$|z(t)|\le\Lambda/2$ for $t\le T$.

\medskip\noindent\emph{Step 1: linear estimates.}
Since $N\ge2$ and $d>0$, we have $A_0,B_0>0$. Thus the roots
$\sigma_1^{(0)},\sigma_2^{(0)}$ of $P_0$ have negative real parts, and
$P_0(\mu)>0$ for every $\mu>0$. Let $G_0$ solve
$P_0(\partial_\tau)G_0=0$ with $G_0(0)=0$ and $G_0'(0)=1$.
Explicitly, for $\tau\ge0$,
\[
 G_0(\tau)=
 \begin{cases}
 \displaystyle\frac{e^{\sigma_1^{(0)}\tau}-e^{\sigma_2^{(0)}\tau}}
 {\sigma_1^{(0)}-\sigma_2^{(0)}},&\sigma_1^{(0)}\ne\sigma_2^{(0)},\\[6pt]
 \tau e^{\sigma_1^{(0)}\tau},&\sigma_1^{(0)}=\sigma_2^{(0)}.
 \end{cases}
\]
The first expression is real also for a complex-conjugate pair.
In all cases, there are constants $C,c>0$ such that
\begin{equation}\label{eq:radial-green-bound}
 |G_0(\tau)|+|G_0'(\tau)|+|G_0''(\tau)|
 \le C(1+\tau)e^{-c\tau}\qquad(\tau\ge0).
\end{equation}
For $g\in C(( -\infty,T])$ with $g(t)=O(e^{\nu t})$, $\nu>0$, define
\begin{equation}\label{QD}
 y_g(t)=\int_{-\infty}^{t}G_0(t-s)g(s)\,\mathrm{d}s.
\end{equation}
The bound \eqref{eq:radial-green-bound} gives absolute convergence and
\[
 |y_g(t)|\le CC_g e^{\nu t}
 \int_0^\infty(1+\tau)e^{-(c+\nu)\tau}\,\mathrm{d}\tau
 =O(e^{\nu t}).
\]
It also permits differentiation under the integral sign. Since
$G_0(0)=0$ and $G_0'(0)=1$,
\[
 y_g'(t)=\int_{-\infty}^tG_0'(t-s)g(s)\,\mathrm{d}s,
 \qquad
 y_g''(t)=g(t)+\int_{-\infty}^tG_0''(t-s)g(s)\,\mathrm{d}s,
\]
so $P_0(\partial_t)y_g=g$. If $y$ is any other solution tending to zero
as $t\to-\infty$, then $y-y_g$ is a homogeneous solution tending to zero.
Because $\operatorname{Re}\sigma_j^{(0)}<0$, the homogeneous solutions
$e^{\sigma_1^{(0)}t},e^{\sigma_2^{(0)}t}$ (or $e^{\sigma_1^{(0)}t},te^{\sigma_1^{(0)}t}$ for a
repeated root) admit no nonzero linear combination with this limit.
Hence $y=y_g$. We have therefore proved that \eqref{QD} is the unique
solution tending to zero at $-\infty$, and that
\[
 P_0(\partial_t)y=O(e^{\nu t}),\quad y(t)\to0
 \quad\Longrightarrow\quad y(t)=O(e^{\nu t})
 \qquad(\nu>0).
\]
In particular, uniqueness gives
\begin{equation}\label{eq:green-exponential}
 \int_{-\infty}^{t}G_0(t-s)e^{\nu s}\,\mathrm{d}s
 =\frac{e^{\nu t}}{P_0(\nu)}\qquad(\nu>0).
\end{equation}
More generally, subtracting $bP_0(\mu)^{-1}e^{\mu t}$ and using the
same estimate yields
\[
 \begin{aligned}
 P_0(\partial_t)y&=be^{\mu t}+O(e^{\nu t}),\quad
 y(t)\to0,\quad 0<\mu<\nu\\
 &\Longrightarrow\quad
 y(t)=\frac{b}{P_0(\mu)}e^{\mu t}+O(e^{\nu t}).
 \end{aligned}
\]
Here and below, all limits are taken as $t\to-\infty$.

\medskip\noindent\emph{Step 2: the first term.}
If $z=O(e^{\eta t})$ for some $\eta>0$, then \eqref{eq:radial-ode} and
$f(z)=O(z^2)$ give
\[
 P_0(\partial_t)z=O\bigl(e^{\min\{2\eta,\rho\}t}\bigr).
\]
Step~1 therefore gives $z=O(e^{\min\{2\eta,\rho\}t})$.
Starting with $\eta_0=\min\{\varepsilon,\rho\}>0$ and iterating
$\eta_{j+1}=\min\{2\eta_j,\rho\}$, we reach $\eta_j=\rho$ after
finitely many steps. Thus
\begin{equation}\label{QDD}
 z(t)=O(e^{\rho t})\qquad\text{as }t\to-\infty.
\end{equation}
Consequently, $f(z)=O(e^{2\rho t})$, and the extraction estimate in
Step~1 gives
\[
 z(t)=d_1e^{\rho t}+O(e^{2\rho t}),\qquad
 d_1=\frac{F}{P_0(\rho)}.
\]
Since $d+\rho=2$, direct calculation gives
$P_0(\rho)=2N+p\,d(N-2+d)>0$; in particular, $d_1\ne0$ and $d_1$ is
independent of $\lambda$.

\medskip\noindent\emph{Step 3: coefficients and remainder estimates.}
For every fixed integer $M\ge2$, Taylor's formula gives, uniformly for
$|s|\le\Lambda/2$,
\begin{equation}\label{eq:finite-taylor-radial}
 f(s)=\sum_{m=2}^M a_ms^m+O(|s|^{M+1}),
\end{equation}
with $a_m$ as in the statement. For $n\ge2$, set
\[
 b_n=\sum_{m=2}^n a_m
 \sum_{\substack{j_1+\cdots+j_m=n\\j_1,\ldots,j_m\ge1}}
 d_{j_1}\cdots d_{j_m}.
\]
Each index $j_i$ is at most $n-1$, so $b_n$ involves only
$d_1,\ldots,d_{n-1}$. Since $P_0(n\rho)>0$, the recurrence
\eqref{eq:radial-recurrence}, namely $d_n=b_n/P_0(n\rho)$ for $n\ge2$,
is well defined and uniquely determines the coefficients.
For $K\ge1$, put
\[
 S_K(t)=\sum_{n=1}^K d_ne^{n\rho t},\qquad \mathcal R_K(t)=z(t)-S_K(t).
\]
We prove by induction that
\begin{equation}\label{eq:radial-induction}
 \mathcal R_K(t)=O(e^{(K+1)\rho t}).
\end{equation}
The case $K=1$ was proved in Step~2. Suppose that
\eqref{eq:radial-induction} holds for some $K\ge1$. Since
$z,S_K=O(e^{\rho t})$, both lie in $[-\Lambda/2,\Lambda/2]$ for
sufficiently negative $t$. Using $f'(s)=O(s)$, we obtain
\[
 \begin{aligned}
 f(z)-f(S_K)
 &=\mathcal R_K\int_0^1 f'(S_K+s\mathcal R_K)\,\mathrm{d}s\\
 &=O(e^{\rho t})\mathcal R_K=O(e^{(K+2)\rho t}).
 \end{aligned}
\]
On the other hand, applying \eqref{eq:finite-taylor-radial} with
$M=K+1$ to $S_K$ and expanding its finite powers gives
\[
 f(S_K)=\sum_{n=2}^{K+1}b_ne^{n\rho t}+O(e^{(K+2)\rho t}).
\]
Indeed, for $2\le n\le K+1$ the coefficient of $e^{n\rho t}$ is
exactly the sum defining $b_n$, since every $j_i\le n-1\le K$.
All omitted products and the Taylor remainder are
$O(e^{(K+2)\rho t})$. The recurrence gives
\[
 P_0(\partial_t)S_K=Fe^{\rho t}+\sum_{n=2}^K b_ne^{n\rho t},
\]
where the sum is empty if $K=1$. Subtracting this identity from
\eqref{eq:radial-ode} and using the preceding estimates yields
\[
 P_0(\partial_t)\mathcal R_K
 =b_{K+1}e^{(K+1)\rho t}+O(e^{(K+2)\rho t}).
\]
Since $\mathcal R_K(t)\to0$ as $t\to-\infty$, Step~1 implies
\[
 \mathcal R_K(t)=\frac{b_{K+1}}{P_0((K+1)\rho)}e^{(K+1)\rho t}
       +O(e^{(K+2)\rho t})
 =d_{K+1}e^{(K+1)\rho t}+O(e^{(K+2)\rho t}).
\]
This is \eqref{eq:radial-induction} with $K$ replaced by $K+1$,
completing the induction.
Finally, $d_1$ is linear in $F$, and each product in $b_n$ has
$j_1+\cdots+j_m=n$. Induction in \eqref{eq:radial-recurrence} therefore
gives $d_n=F^nD_n(\lambda,p,\alpha,N)$.
\end{proof}

\section{Asymptotic behavior of asymptotically radial rupture solutions: the nonradial case}

We retain the notation of Section~1. The proof of
Proposition~\ref{prop:expansion} starts with the independent linear estimate
in Lemma~\ref{lem:linear-extraction}, uses it to obtain the initial decay
in Proposition~\ref{prop:decay}, and then extracts successive terms of the
expansion. Fix $T_R<\min\{0,\ln R\}$ and set
$t=\ln|x|$, $\theta=x/|x|$, and
$z(t,\theta)=|x|^{-d}u(x)-\Lambda$. For a positive asymptotically radial
rupture solution, $z=O(e^{\varepsilon t})$ as $t\to-\infty$, and \eqref{3}
becomes
\[
\mathcal Lz=f(z)+q_F(t),\qquad
q_F(t)=
\begin{cases}
Fe^{\rho t},&F\ne0,\\
0,&F=0,
\end{cases}
\]
where
\[
f(z):=\frac{\lambda}{(\Lambda+z)^p}
-\frac{\lambda}{\Lambda^p}
+\frac{\lambda p}{\Lambda^{p+1}}z
=O(z^2)\quad\text{as }z\to0,
\]
so $f$ is analytic near zero and $f(0)=f'(0)=0$. Here, with $A_0,B_0$ as in \eqref{eq:AB},
\begin{equation}\label{eq:L}
\mathcal L=\partial_t^2+A_0\partial_t+\Delta_{S^{N-1}}+B_0.
\end{equation}
The modal operators are
\begin{equation}\label{wqwqq}
\mathcal L_k=\frac{\mathrm d^2}{\mathrm dt^2}+A_0\frac{\mathrm d}{\mathrm dt}+B_0-\lambda_k,
\qquad \lambda_k=k(N-2+k).
\end{equation}
Let $\{Q_1^k,\ldots,Q_{m_k}^k\}$ be an orthonormal basis of the degree-$k$
eigenspace, where $m_0=1$ and
\[
m_k=\frac{(N-2+2k)(N-3+k)!}{k!(N-2)!}\quad(k\ge1).
\]
The indicial polynomial of $\mathcal L_k$ is $P_k$ in
\eqref{eq:modal-indicial-polynomial}, with roots
$\sigma_1^{(k)},\sigma_2^{(k)}$ and threshold $k_*$ given by
\eqref{eq:modal-roots}--\eqref{eq:kstar}.
For $k\ge k_*$, the roots are real and satisfy
$\sigma_1^{(k)}>0>\sigma_2^{(k)}$; for $k<k_*$, both have nonpositive
real part.

The notation below describes the exponents and angular degrees in
Proposition~\ref{prop:expansion}. Decaying homogeneous modes and the
pressure term (when $F\ne0$) supply the generators; nonlinear products
add their exponents. The initial decay proof uses only the spectral data
and $\mu_1$, not the remaining expansion notation.

\medskip\noindent\textbf{The exponent set.}
The positive indicial roots $\gamma_k=\sigma_1^{(k)}$, $k\ge k_*$,
increase strictly to infinity. The generating set $\mathcal G_F$ in
\eqref{eq:GF} also includes $\rho>0$ when $F\ne0$. Thus
$\min\mathcal G_F>0$ and $\mathcal G_F$ has no finite accumulation point.
Let $\mathbb N_0^{(\mathcal G_F)}$ denote the set of finitely supported
families $\mathbf n=(n_\gamma)_{\gamma\in\mathcal G_F}$ of nonnegative
integers, and write $|\mathbf n|=\sum_\gamma n_\gamma$ and
$\langle\mathbf n\rangle=\sum_\gamma n_\gamma\,\gamma$.
In this notation the semigroup $\mathcal I$ of \eqref{eq:I-sets} is
$\mathcal I=\{\langle\mathbf n\rangle:|\mathbf n|\ge1\}$, and we split it as
\begin{equation}\label{eq:I-split}
\mathcal I_\rho=\mathcal G_F=\{\langle\mathbf n\rangle:|\mathbf n|=1\},\qquad
\mathcal I_{\widetilde\rho}=\{\langle\mathbf n\rangle:|\mathbf n|\ge2\},
\end{equation}
so that $\mathcal I=\mathcal I_\rho\cup\mathcal I_{\widetilde\rho}$, and
$\mathcal I_{\widetilde\rho}$ consists of sums of at least two generators,
with repetitions allowed. If $\langle\mathbf n\rangle\le M$, then
$n_\gamma=0$ for $\gamma>M$ and $n_\gamma\le M/\min\mathcal G_F$.
Since $\mathcal G_F\cap(0,M]$ is finite, only finitely many such
$\mathbf n$ exist. Hence each bounded interval meets $\mathcal I$ in a
finite set, and $\mathcal I$ and $\mathcal I_{\widetilde\rho}$ can be
enumerated as strictly increasing sequences
\[
\mu_1<\mu_2<\cdots\to\infty,\qquad
\widetilde\rho_1<\widetilde\rho_2<\cdots\to\infty,
\]
the first of which is the sequence \eqref{eq:mu-seq}, with
$\mu_1=\min\mathcal G_F$ and $\widetilde\rho_1=2\mu_1$; explicitly,
$\mu_1=\gamma_{k_*}$ if $F=0$ and
$\mu_1=\min\{\gamma_{k_*},\rho\}$ if $F\ne0$.
The intersection $\mathcal I_\rho\cap\mathcal I_{\widetilde\rho}$ need
not be empty: a generator may itself be a sum of two or more generators.
An overlap of exponents alone does not imply modal resonance, which
requires a nonzero forcing component of degree $\ell$ at an exponent
$\nu$ with $P_\ell(\nu)=0$.

\smallskip\noindent\textbf{Angular-degree bounds.}
To each $\gamma\in\mathcal G_F$ we attach its \emph{angular degree}
\begin{equation}\label{eq:angular-degree}
 \mathfrak d(\gamma)=
 \begin{cases}
 k, & \text{if } \gamma=\gamma_k \text{ for some } k\ge k_*,\\[2pt]
 0, & \text{if } F\ne0,\ \gamma=\rho \text{ and } \rho\notin\{\gamma_k:k\ge k_*\}.
 \end{cases}
\end{equation}
Thus $\mathfrak d(\gamma_k)=k$ and the pressure generator has degree $0$
unless it coincides with some $\gamma_j$, in which case the common
generator is assigned degree $j$. This records the largest
spherical-harmonic degree contributed by a single generator.
For $\nu\in\mathcal I_{\widetilde\rho}$ set
\begin{equation}\label{eq:Mtilde}
\widetilde M(\nu)=\max\Big\{\sum_{\gamma\in\mathcal G_F}n_\gamma\,\mathfrak d(\gamma)\ :\ \mathbf n\in\mathbb N_0^{(\mathcal G_F)},\ |\mathbf n|\ge2,\ \langle\mathbf n\rangle=\nu\Big\},
\end{equation}
a well-defined nonnegative integer, since the maximum is taken over a
finite nonempty set. By Lemma~\ref{lem:harmonics} below and induction on
the number of factors, $\widetilde M(\nu)$ bounds the spherical-harmonic
degree of every product $\prod_\gamma\Theta_\gamma^{\,n_\gamma}$ in which
each $\Theta_\gamma$ is a finite sum of spherical harmonics of degree at
most $\mathfrak d(\gamma)$ and $|\mathbf n|\ge2$,
$\langle\mathbf n\rangle=\nu$.
For $\mu\in\mathcal I$, also define
\[
 M^*(\mu)=
 \max\Big\{\sum_{\gamma\in\mathcal G_F}n_\gamma\,\mathfrak d(\gamma):
 \mathbf n\in\mathbb N_0^{(\mathcal G_F)},\ |\mathbf n|\ge1,\ \langle\mathbf n\rangle=\mu\Big\}.
\]
By local finiteness, the maximum is finite. It equals $\mathfrak d(\mu)$ if
$\mu\notin\mathcal I_{\widetilde\rho}$, $\widetilde M(\mu)$ if
$\mu\notin\mathcal G_F$, and their maximum if $\mu$ belongs to both sets.

\begin{prop}\label{prop:expansion}
Assume \eqref{H}, and let $d$ and $\Lambda$ be given by \eqref{qcc}.
Let $u\in C^2(B_R\setminus\{0\})$ be a positive rupture solution of
\eqref{QQ} such that $u(x)=\Lambda|x|^d(1+O(|x|^\varepsilon))$ as
$|x|\to0^+$ for some $\varepsilon>0$.
Set $z(t,\theta)=|x|^{-d}u(x)-\Lambda$, with $t=\ln|x|$ and
$\theta=x/|x|$, and let $\{\mu_j\}_{j\ge1}$ be the sequence in
\eqref{eq:mu-seq}. There exist nonnegative integers $M_j$ and coefficient
functions $C_{ji}$, $j\ge1$, $0\le i\le M_j$, independent of the
truncation order, such that, for every positive integer $n$, uniformly
in $\theta\in S^{N-1}$,
\[
z(t,\theta)=\sum_{j=1}^n\sum_{i=0}^{M_j}C_{ji}(\theta)t^ie^{\mu_jt}
+O\left((1+|t|)^{M_{n+1}}e^{\mu_{n+1}t}\right),
\qquad t\to-\infty.
\]
Each $C_{ji}$ is a finite sum of spherical harmonics of degree at most
$M^*(\mu_j)$; equivalently,
$C_{ji}(\theta)=\sum_{k=0}^{M^*(\mu_j)}\sum_{q=1}^{m_k}a_{jikq}Q_q^k(\theta)$
in the fixed spherical-harmonic basis, with constant coefficients $a_{jikq}$.
\end{prop}

\begin{prop}\label{prop:decay}
Assume \eqref{H}, and let $d$ and $\Lambda$ be given by \eqref{qcc}.
Let $u\in C^2(B_R\setminus\{0\})$ be a positive rupture solution of
\eqref{QQ} satisfying $u(x)=\Lambda|x|^d(1+O(|x|^\varepsilon))$ as
$|x|\to0^+$ for some $\varepsilon>0$.
Define $z(t,\theta)=|x|^{-d}u(x)-\Lambda$, with $t=\ln|x|$ and
$\theta=x/|x|$, and let $\mu_1$ be given by \eqref{eq:mu-seq}.
For the fixed $T_R<\min\{0,\ln R\}$, there exists a constant $C>0$ such that
\[
 \max_{\theta\in S^{N-1}}|z(t,\theta)|\le Ce^{\mu_1t},
 \qquad t\le T_R.
\]
\end{prop}

\medskip\noindent\textbf{A common linear estimate.}
As $t\to-\infty$, homogeneous modes with $0<\gamma_k<\Gamma_*$ decay
more slowly than the forcing bound $(1+|t|)^{m_*}e^{\Gamma_*t}$.
The next lemma extracts these modes and bounds the remainder by the same
weight, with a possible extra factor $1+|t|$ at resonance
$\Gamma_*=\gamma_k$.

\begin{lem}[Extraction of decaying homogeneous modes]\label{lem:linear-extraction}
Let $T<0$, $\Gamma_*>0$, and $m_*\in\mathbb N_0$. Suppose that
$g$ is continuous on $(-\infty,T]\times S^{N-1}$ and
\[
 |g(t,\theta)|\le C(1+|t|)^{m_*} e^{\Gamma_* t},
 \qquad (t,\theta)\in(-\infty,T]\times S^{N-1}.
\]
If $v\in C^2(( -\infty,T]\times S^{N-1})$ satisfies
\[
 \mathcal Lv=g,
 \qquad
 v(t,\theta)\longrightarrow0\quad\text{as }t\to-\infty
 \quad\text{uniformly in }\theta,
\]
then there are constants $c_j^k$, indexed by the finitely many pairs
\[
 k\ge k_*,\qquad \gamma_k<\Gamma_*,\qquad 1\le j\le m_k,
\]
such that
\begin{equation}\label{eq:linear-extraction}
 \left|v(t,\theta)-
 \sum_{\substack{k\ge k_*:\ \gamma_k<\Gamma_*}}
 \sum_{j=1}^{m_k}c_j^k e^{\gamma_kt}Q_j^k(\theta)\right|
 \le C'(1+|t|)^{m_*+\chi_{\Gamma_*}}e^{\Gamma_* t}
\end{equation}
for all $(t,\theta)\in(-\infty,T]\times S^{N-1}$, where
\[
 \chi_{\Gamma_*}=
 \begin{cases}
 1,&\Gamma_*=\gamma_k\text{ for some }k\ge k_*,\\
 0,&\Gamma_*\notin\{\gamma_k:k\ge k_*\}.
 \end{cases}
\]
The sum in \eqref{eq:linear-extraction} is interpreted as zero when
$0<\Gamma_*\le\gamma_{k_*}$. The constant $C'$ is independent of $t$ and
$\theta$, but may depend on $N,A_0,B_0,T,\Gamma_*,m_*,C$, and $v$.
\end{lem}

\begin{proof}
We first estimate the individual spherical-harmonic modes, then control
the high-mode tail in $L^2$, and finally apply local elliptic estimates
to obtain \eqref{eq:linear-extraction}. Put
\[
 \Xi_{\Gamma_*}(t)=(1+|t|)^{m_*}e^{\Gamma_*t},\qquad t\le T.
\]

\smallskip\noindent
\emph{Step 1. Fixed-mode estimates.}
Let $v_j^k$ and $g_j^k$ be the coefficients of $v$ and $g$ in the
orthonormal basis $\{Q_j^k\}$. Then
$P_k(\partial_t)v_j^k=g_j^k$ and $|g_j^k(t)|\le C\Xi_{\Gamma_*}(t)$.
For $0\le k<k_*$, define $G_k$ on $[0,\infty)$ by
\[
 P_k(\partial_s)G_k=0,\qquad G_k(0)=0,\qquad G_k'(0)=1.
\]
The roots of $P_k$ have nonpositive real parts, so
$|G_k(s)|+|G_k'(s)|+|G_k''(s)|\le C_k(1+s)$ for $s\ge0$.
This bound also covers repeated, complex-conjugate, and zero roots.
Consequently,
\[
 v_j^k(t)=\int_{-\infty}^tG_k(t-s)g_j^k(s)\,\mathrm{d}s.
\]
Indeed, the integral converges, solves the modal equation, and tends to
zero at $-\infty$. Its difference from $v_j^k$ is a homogeneous solution
with the same limit; since no nonzero homogeneous solution tends to zero
at $-\infty$, this difference vanishes.
Changing variables to $\tau=t-s\ge0$ and using
$1+|t-\tau|\le(1+|t|)(1+\tau)$, we obtain
\begin{equation}\label{eq:low-extraction-bound}
 |v_j^k(t)|\le C_k\Xi_{\Gamma_*}(t)
 \int_0^\infty(1+\tau)^{m_*+1}e^{-\Gamma_*\tau}\,\mathrm{d}\tau
 \le C_k\Xi_{\Gamma_*}(t),\qquad 0\le k<k_*.
\end{equation}

For $k\ge k_*$, the roots satisfy
$\gamma_k>0>\sigma_2^{(k)}$ and
$\gamma_k-\sigma_2^{(k)}=\sqrt{\Delta_k}$. Variation of constants,
together with $v_j^k(t)\to0$ as $t\to-\infty$, gives constants $a_j^k$
such that
\[
\begin{aligned}
 v_j^k(t)={}&a_j^k e^{\gamma_kt}
 -\dfrac{1}{\sqrt{\Delta_k}}\int_t^T e^{\gamma_k(t-s)}g_j^k(s)\,\mathrm{d}s\\
 &-\dfrac{1}{\sqrt{\Delta_k}}\int_{-\infty}^t
 e^{\sigma_2^{(k)}(t-s)}g_j^k(s)\,\mathrm{d}s.
\end{aligned}
\]
The assumed decay of $g_j^k$ makes both integrals converge and tend to
zero as $t\to-\infty$; thus the additional homogeneous term
$e^{\sigma_2^{(k)}t}$ has zero coefficient. The negative-root integral
satisfies
\[
 \left|\int_{-\infty}^t e^{\sigma_2^{(k)}(t-s)}g_j^k(s)\,\mathrm{d}s\right|
 \le C\Xi_{\Gamma_*}(t)\int_0^\infty
 (1+s)^{m_*}e^{-(\Gamma_*-\sigma_2^{(k)})s}\,\mathrm{d}s
 \le C_k\Xi_{\Gamma_*}(t).
\]
If $\gamma_k<\Gamma_*$, the integral defining
\[
 c_j^k=a_j^k-\dfrac{1}{\sqrt{\Delta_k}}
 \int_{-\infty}^T e^{-\gamma_ks}g_j^k(s)\,\mathrm{d}s
\]
converges, and subtraction gives
\[
 v_j^k(t)-c_j^ke^{\gamma_kt}
 =\dfrac{1}{\sqrt{\Delta_k}}\int_{-\infty}^t
 \left(e^{\gamma_k(t-s)}-e^{\sigma_2^{(k)}(t-s)}\right)
 g_j^k(s)\,\mathrm{d}s
 =O(\Xi_{\Gamma_*}(t)).
\]
The positive-root integral has the same bound, with
$\Gamma_*-\gamma_k>0$ in place of $\Gamma_*-\sigma_2^{(k)}$.
If $\gamma_k\ge\Gamma_*$, then
$(1+|s|)^{m_*}\le(1+|t|)^{m_*}$ for $t\le s\le T$, so
\[
\begin{aligned}
 \left|\int_t^T e^{\gamma_k(t-s)}g_j^k(s)\,\mathrm{d}s\right|
 &\le C\Xi_{\Gamma_*}(t)\int_0^{T-t}e^{-(\gamma_k-\Gamma_*)s}\,\mathrm{d}s\\
 &\le
 \begin{cases}
 C(1+|t|)\Xi_{\Gamma_*}(t),&\gamma_k=\Gamma_*,\\[2mm]
 \dfrac{C}{\gamma_k-\Gamma_*}\Xi_{\Gamma_*}(t),&\gamma_k>\Gamma_*.
 \end{cases}
\end{aligned}
\]
The homogeneous term is also $O(\Xi_{\Gamma_*}(t))$ when
$\gamma_k\ge\Gamma_*$. Consequently,
\begin{equation}\label{eq:fixed-mode-extraction}
\begin{cases}
 |v_j^k(t)-c_j^ke^{\gamma_kt}|
 \le C_k(1+|t|)^{m_*} e^{\Gamma_* t},&\gamma_k<\Gamma_*,\\[2mm]
 |v_j^k(t)|
 \le C_k(1+|t|)^{m_*+1} e^{\Gamma_* t},&\gamma_k=\Gamma_*,\\[2mm]
 |v_j^k(t)|
 \le C_k(1+|t|)^{m_*} e^{\Gamma_* t},&\gamma_k>\Gamma_*.
\end{cases}
\end{equation}

\smallskip\noindent
\emph{Step 2. An $L^2$ estimate for the high-mode tail.}
The constants in \eqref{eq:fixed-mode-extraction} depend on the mode, so
we estimate the infinite tail directly. Choose $\mathcal K\ge k_*$ such
that $\gamma_{\mathcal K+1}>\Gamma_*$. This cutoff is used only for the
tail estimate; the terms to be extracted are still precisely those with
$\gamma_k<\Gamma_*$. Let $\Pi_{>\mathcal K}$ be the orthogonal projection
onto the spherical harmonics of degrees greater than $\mathcal K$, and set
\[
 \widehat v=\Pi_{>\mathcal K}v,\qquad
 \widehat g=\Pi_{>\mathcal K}g.
\]
Below, $\|\cdot\|_2$ and $\|\cdot\|_\infty$ denote norms on $S^{N-1}$.
Since this projection commutes with $\mathcal L$ and subtracts only a fixed
finite-dimensional projection,
\[
 \mathcal L\widehat v=\widehat g,\qquad
 \|\widehat g(t,\cdot)\|_2+\|\widehat g(t,\cdot)\|_\infty
 \le C\Xi_{\Gamma_*}(t).
\]
For $k>\mathcal K$ we have
$\gamma_k\ge\gamma_{\mathcal K+1}$,
$\sigma_2^{(k)}\le\sigma_2^{(\mathcal K+1)}<0$, and
$\Delta_k\ge\Delta_{\mathcal K+1}>0$.
Evaluating the representation in Step~1 at $t=T$ and using Parseval's
identity and Minkowski's inequality gives
\[
 \left(\sum_{k>\mathcal K}\sum_{j=1}^{m_k}
 |a_j^ke^{\gamma_kT}|^2\right)^{1/2}
 \le \|\widehat v(T,\cdot)\|_2
 +C\int_{-\infty}^{T}
 e^{\sigma_2^{(\mathcal K+1)}(T-s)}\Xi_{\Gamma_*}(s)\,\mathrm{d}s<\infty.
\]
These inequalities extend from finite sums by passage to the limit.
For $t\le T$, the coefficient bound and
$e^{\gamma_k(t-T)}\le e^{\gamma_{\mathcal K+1}(t-T)}$
control the homogeneous tail; Minkowski's inequality for the integral
terms gives
\begin{equation}\label{eq:tail-integral-bound}
\begin{aligned}
 \|\widehat v(t,\cdot)\|_2
 \le{}&Ce^{\gamma_{\mathcal K+1}(t-T)}
 +C\int_t^T e^{\gamma_{\mathcal K+1}(t-s)}\Xi_{\Gamma_*}(s)\,\mathrm{d}s\\
 &+C\int_{-\infty}^t
 e^{\sigma_2^{(\mathcal K+1)}(t-s)}\Xi_{\Gamma_*}(s)\,\mathrm{d}s.
\end{aligned}
\end{equation}
By the weighted integral estimates in Step~1, the two integrals are
bounded, respectively, by
\[
 \dfrac{\Xi_{\Gamma_*}(t)}{\gamma_{\mathcal K+1}-\Gamma_*}
 \quad\text{and}\quad
 \Xi_{\Gamma_*}(t)\int_0^\infty(1+s)^{m_*}
 e^{-(\Gamma_*-\sigma_2^{(\mathcal K+1)})s}\,\mathrm{d}s.
\]
Both are $O(\Xi_{\Gamma_*}(t))$. Since
$\gamma_{\mathcal K+1}>\Gamma_*$, the homogeneous term has this bound as
well. Hence
\begin{equation}\label{eq:tail-L2-extraction}
 \|\widehat v(t,\cdot)\|_2
 \le C\Xi_{\Gamma_*}(t)
 =C(1+|t|)^{m_*}e^{\Gamma_*t},\qquad t\le T.
\end{equation}

\smallskip\noindent
\emph{Step 3. Uniform estimate and conclusion.}
The standard interior $L^\infty$ estimate for the uniformly elliptic
equation $\mathcal L\widehat v=\widehat g$
gives, for $t\le T-2$,
\[
 \|\widehat v(t,\cdot)\|_\infty
 \le C\left[
 \left(\int_{t-1}^{t+1}\|\widehat v(s,\cdot)\|_2^2\,\mathrm{d}s\right)^{1/2}
 +\sup_{|s-t|\le1}\|\widehat g(s,\cdot)\|_\infty\right].
\]
The constant $C$ is independent of $t$, since the coefficients of
$\mathcal L$ are independent of $t$ and $S^{N-1}$ admits a finite
covering by fixed coordinate charts.
By \eqref{eq:tail-L2-extraction} and the pointwise bound for $\widehat g$,
both terms on the right are bounded by $C\Xi_{\Gamma_*}(t)$, since the
weight is uniformly comparable on unit intervals. Thus
\begin{equation}\label{eq:tail-Linfty-extraction}
 |\widehat v(t,\theta)|\le C(1+|t|)^{m_*}e^{\Gamma_*t},
 \qquad t\le T-2.
\end{equation}
Increasing $C$ covers the remaining compact interval.
The modes of degree at most $\mathcal K$ form a finite sum, so
\eqref{eq:low-extraction-bound} and \eqref{eq:fixed-mode-extraction}
can now be summed over them. Since $\gamma_{\mathcal K+1}>\Gamma_*$,
this finite part contains every extracted or resonant mode. Adding
\eqref{eq:tail-Linfty-extraction} proves \eqref{eq:linear-extraction};
the possible extra factor $1+|t|$ comes only from the case
$\gamma_k=\Gamma_*$ in \eqref{eq:fixed-mode-extraction}.
\end{proof}

\medskip\noindent\textbf{The initial decay estimate.}
\begin{proof}[Proof of Proposition~\ref{prop:decay}]
All estimates are uniform in $\theta$ as $t\to-\infty$.
Remove the explicit pressure forcing by setting
\[
 \eta_\star(t)=
 \begin{cases}
 \dfrac{F}{P_0(\rho)}e^{\rho t},&F\ne0,\\[4pt]
 0,&F=0,
 \end{cases}
 \qquad v=z-\eta_\star.
\]
For $F\ne0$, $\rho>0$ and $P_0(\rho)=\rho^2+A_0\rho+B_0>0$.
Thus $\mathcal L\eta_\star=q_F$ and
$\eta_\star=O(e^{\mu_1t})$, since $\rho\ge\mu_1$ when $F\ne0$.
In both cases,
\[
 \mathcal Lv=f(z),\qquad v\longrightarrow0\quad(t\to-\infty).
\]
The right-hand side remains $f(z)$, not $f(v)$.

\smallskip\noindent
\emph{Step 1. Improve the decay below the first exponent.}
Choose $\widehat\mu=3\mu_1/4$, so that
$\widehat\mu<\mu_1<2\widehat\mu$.
This keeps the iteration below the first positive indicial root, while
the strict gap $2\widehat\mu-\mu_1>0$ absorbs any polynomial loss in the
final step. Set
\[
 \kappa_0=\min\{\varepsilon,\widehat\mu\},\qquad
 \kappa_{n+1}=\min\{2\kappa_n,\widehat\mu\}.
\]
Initially $z=O(e^{\kappa_0t})$. If $z=O(e^{\kappa_nt})$, then
\[
 \mathcal Lv=f(z)=O(e^{2\kappa_nt})
 =O(e^{\kappa_{n+1}t}).
\]
Apply Lemma~\ref{lem:linear-extraction} with
$\Gamma_*=\kappa_{n+1}$ and $m_*=0$.
Since $0<\kappa_{n+1}\le\widehat\mu<\mu_1\le\gamma_{k_*}$,
there is neither an extracted mode nor a resonant factor, and hence
$v=O(e^{\kappa_{n+1}t})$. Adding $\eta_\star$ gives the same bound
for $z$. Therefore $\kappa_n=\min\{2^n\kappa_0,\widehat\mu\}$,
and after finitely many steps
\[
 z=O(e^{\widehat\mu t}).
\]

\smallskip\noindent
\emph{Step 2. Reach the exponent $\mu_1$ without polynomial loss.}
Now $\mathcal Lv=f(z)=O(e^{2\widehat\mu t})$.
A further application of Lemma~\ref{lem:linear-extraction}, with
$\Gamma_*=2\widehat\mu$ and $m_*=0$, gives
\[
 v=\sum_{\substack{k\ge k_*:\ \gamma_k<2\widehat\mu}}
 \sum_{j=1}^{m_k}c_j^k e^{\gamma_kt}Q_j^k
 +O\!\left((1+|t|)^{\chi_{2\widehat\mu}}
 e^{2\widehat\mu t}\right).
\]
Every exponent in this finite sum satisfies
$\gamma_k\ge\gamma_{k_*}\ge\mu_1$.
Moreover, $2\widehat\mu=3\mu_1/2>\mu_1$, so even in the resonant case
\[
 (1+|t|)^{\chi_{2\widehat\mu}}e^{2\widehat\mu t}
 =O(e^{\mu_1t}).
\]
Hence $v=O(e^{\mu_1t})$, and adding back $\eta_\star$ yields
$z=O(e^{\mu_1t})$. Continuity extends the estimate to all $t\le T_R$.
\end{proof}

\medskip\noindent\textbf{The higher-order expansion.}
We use the following elementary product rule
(see \cite[Lemma~21]{8z}, \cite{5z}).
\begin{lem}\label{lem:harmonics}
If $Q_k$ and $Q_l$ are spherical harmonics of degrees $k$ and $l$ respectively, then
\[
 Q_kQ_l=\sum_{i=0}^{k+l}Z_i,
\]
where each $Z_i$ is either zero or a spherical harmonic of degree $i$.
\end{lem}

\begin{proof}[Proof of Proposition~\ref{prop:expansion}]
Proposition~\ref{prop:decay} gives $z=O(e^{\mu_1t})$. We construct the
expansion successively along $\{\mu_j\}$, controlling both the remainder
and the residual of each partial sum. All estimates below hold as
$t\to-\infty$, uniformly in $\theta\in S^{N-1}$. The powers of the
polynomial weights are finite nonnegative integers and may depend on the
expansion order. We first show how to extract one exponent, then iterate
this step. At each finite stage, we work sufficiently far to the left
that the arguments of $f$ lie in its analytic neighborhood of zero.

For $\nu>0$, let $\mathscr P_\nu$ denote the space of functions
\[
 e^{\nu t}\sum_{q=0}^{\mathcal J}t^q\Theta_q(\theta),
\]
where $\mathcal J<\infty$ and each $\Theta_q$ is a finite sum of spherical
harmonics.

\smallskip\noindent
\emph{Step 1. One-step extraction.}
Let $\nu<\nu^+$ be consecutive elements of $\mathcal I$. Suppose that
$\mathcal R\in C^2(( -\infty,T_R]\times S^{N-1})$ and, for some
$A,B\in\mathbb N_0$ and $\mathscr S_\nu\in\mathscr P_\nu$,
\begin{equation}\label{eq:linear-extraction-hyp}
 \mathcal R=O\!\left((1+|t|)^A e^{\nu t}\right),\qquad
 \mathcal L\mathcal R=\mathscr S_\nu+
 O\!\left((1+|t|)^B e^{\nu^+t}\right).
\end{equation}
Then there is a function $\Psi_\nu\in\mathscr P_\nu$ such that
\begin{equation}\label{eq:linear-extraction-concl}
 \mathcal L\Psi_\nu=\mathscr S_\nu,\qquad
 \mathcal R-\Psi_\nu
 =O\!\left((1+|t|)^{A^+} e^{\nu^+t}\right),
\end{equation}
where one may take $A^+=B$ unless $\nu^+$ is a positive indicial root,
in which case one may take $A^+=B+1$.

To prove the claim, decompose $\mathscr S_\nu$ into spherical harmonics.
For a component of degree $\ell$, seek a particular solution in the form
$e^{\nu t}\varpi(t)Q_q^\ell(\theta)$, with $\varpi$ a polynomial. Since
\[
 \mathcal L_\ell\!\left(e^{\nu t}\varpi(t)\right)
 =e^{\nu t}\bigl[\varpi''+P_\ell'(\nu)\varpi'+P_\ell(\nu)\varpi\bigr],
\]
comparison of polynomial coefficients gives a solution of the same degree
as the forcing polynomial when $P_\ell(\nu)\ne0$. If $P_\ell(\nu)=0$,
then $P_\ell'(\nu)=2\nu+A_0>0$, so $\nu$ is a simple positive indicial
root and a polynomial solution has degree at most one higher.
Resonance at $\nu$ concerns the polynomial degree of the extracted term,
whereas resonance at $\nu^+$ concerns the remainder weight $A^+$ in
\eqref{eq:linear-extraction-concl}. Summing these particular solutions gives $\Psi_\nu^{\mathrm{part}}\in\mathscr P_\nu$
with $\mathcal L\Psi_\nu^{\mathrm{part}}=\mathscr S_\nu$. The difference
$\mathcal R-\Psi_\nu^{\mathrm{part}}$ tends to zero uniformly at $-\infty$
and satisfies
\[
 \mathcal L\bigl(\mathcal R-\Psi_\nu^{\mathrm{part}}\bigr)
 =O\!\left((1+|t|)^B e^{\nu^+t}\right).
\]
Apply Lemma~\ref{lem:linear-extraction} to this difference with
$\Gamma_*=\nu^+$ and $m_*=B$. It extracts the homogeneous modes with
positive indicial roots below $\nu^+$. Those with roots below $\nu$ have
zero coefficients: indeed,
$\mathcal R-\Psi_\nu^{\mathrm{part}}=O((1+|t|)^{A'}e^{\nu t})$ for some
finite $A'$. For $\gamma_k<\nu$, project the conclusion of the lemma onto
$Q_j^k$ and multiply by $e^{-\gamma_kt}$. Both this bound and the remainder
at exponent $\nu^+$ tend to zero after multiplication, so $c_j^k=0$.
Every positive indicial root lies in $\mathcal G_F\subset\mathcal I$;
since $\nu$ and $\nu^+$ are consecutive, the only remaining extracted
modes have root exactly $\nu$. Adding them to
$\Psi_\nu^{\mathrm{part}}$ proves \eqref{eq:linear-extraction-concl} and
the claim. In particular, $\Psi_\nu$ may contain a homogeneous term even
when $\mathscr S_\nu=0$.

\smallskip\noindent
\emph{Step 2. Induction over the exponent set.}
We construct $\Psi_{\mu_j}\in\mathscr P_{\mu_j}$ successively. For $n\ge1$,
set
\[
 S_{n-1}=\sum_{j=1}^{n-1}\Psi_{\mu_j},\qquad
 \mathcal R_n=z-S_{n-1}.
\]
The equation $\mathcal Lz=f(z)+q_F$ gives
\[
 \mathcal L\mathcal R_n
 =\bigl[f(S_{n-1})+q_F-\mathcal LS_{n-1}\bigr]
  +\bigl[f(z)-f(S_{n-1})\bigr].
\]
The first bracket is the residual of $S_{n-1}$; the second is controlled
by $\mathcal R_n$ using $f(0)=f'(0)=0$.
We therefore induct on both the remainder and the residual, assuming
that, for some $A_n\in\mathbb N_0$,
\begin{equation}\label{eq:expansion-induction}
\begin{aligned}
 \text{(a)}\quad &\mathcal R_n
 =O\!\left((1+|t|)^{A_n}e^{\mu_nt}\right),\\
 \text{(b)}\quad &f(S_{n-1})+q_F-\mathcal LS_{n-1}
 =O\!\left((1+|t|)^{A_n}e^{\mu_nt}\right).
\end{aligned}
\end{equation}
These hold for $n=1$ with $S_0=0$ and $A_1=0$: condition~(a) is the
initial estimate, and the residual in~(b) is $q_F$, which is either zero
or has exponent $\rho\ge\mu_1$. At each step, $B$ denotes a finite
nonnegative integer that may be increased to control polynomial weights.

Since $f$ is analytic near zero, write $f(y)=\sum_{m=2}^{\infty}a_my^m$
there. Choose $J_n\ge2$ such that $(J_n+1)\mu_1\ge\mu_{n+1}$. Every term
of $S_{n-1}$ has exponent at least $\mu_1$, so Taylor's formula gives
\[
 f(S_{n-1})-\sum_{m=2}^{J_n}a_mS_{n-1}^m
 =O\!\left((1+|t|)^B e^{(J_n+1)\mu_1t}\right).
\]
All exponents in the finite polynomial--exponential sum
\[
 \sum_{m=2}^{J_n}a_mS_{n-1}^m+q_F-\mathcal LS_{n-1}
\]
belong to $\mathcal I$. After equal exponents are grouped together,
condition~\eqref{eq:expansion-induction}(b) excludes every nonzero term
below $\mu_n$: otherwise, dividing by the smallest such exponential
would force a nonzero polynomial in $t$ to tend to zero. Let
$\mathscr S_{\mu_n}\in\mathscr P_{\mu_n}$ be the sum of its terms with
exponent exactly $\mu_n$. Increasing $J_n$ adds only monomials with
exponent at least $(J_n+1)\mu_1\ge\mu_{n+1}>\mu_n$, so
$\mathscr S_{\mu_n}$ is independent of this choice. Since no element of
$\mathcal I$ lies between $\mu_n$ and $\mu_{n+1}$, we obtain
\begin{equation}\label{eq:expansion-leading-residual}
 f(S_{n-1})+q_F-\mathcal LS_{n-1}
 =\mathscr S_{\mu_n}+
 O\!\left((1+|t|)^B e^{\mu_{n+1}t}\right).
\end{equation}

The additive closure of $\mathcal I$ gives
$\mu_1+\mu_n\in\mathcal I$, hence $\mu_1+\mu_n\ge\mu_{n+1}$.
Since $f(0)=f'(0)=0$ and $S_{n-1},\mathcal R_n\to0$ uniformly,
\begin{equation}\label{eq:nonlinear-difference}
\begin{aligned}
 f(z)-f(S_{n-1})
 &=O\!\left((|S_{n-1}|+|\mathcal R_n|)|\mathcal R_n|\right)\\
 &=O\!\left((1+|t|)^B e^{(\mu_1+\mu_n)t}\right).
\end{aligned}
\end{equation}
Substituting \eqref{eq:expansion-leading-residual} and
\eqref{eq:nonlinear-difference} into the equation for $\mathcal R_n$ yields
\[
 \mathcal L\mathcal R_n=\mathscr S_{\mu_n}+
 O\!\left((1+|t|)^B e^{\mu_{n+1}t}\right).
\]
Applying the claim with $\nu=\mu_n$ and $\nu^+=\mu_{n+1}$ gives
$\Psi_{\mu_n}\in\mathscr P_{\mu_n}$ such that
\[
 \mathcal L\Psi_{\mu_n}=\mathscr S_{\mu_n},\qquad
 \mathcal R_{n+1}:=\mathcal R_n-\Psi_{\mu_n}
 =O\!\left((1+|t|)^{A_{n+1}}e^{\mu_{n+1}t}\right)
\]
for some finite $A_{n+1}\in\mathbb N_0$.

Set $S_n=S_{n-1}+\Psi_{\mu_n}$. The residual satisfies the exact identity
\begin{equation}\label{eq:expansion-residual-update}
\begin{aligned}
 f(S_n)+q_F-\mathcal LS_n
 ={}&\bigl[f(S_{n-1})+q_F-\mathcal LS_{n-1}-\mathscr S_{\mu_n}\bigr]\\
 &+\bigl[f(S_n)-f(S_{n-1})\bigr].
\end{aligned}
\end{equation}
The first bracket is controlled by \eqref{eq:expansion-leading-residual}.
The second is
$O((|S_{n-1}|+|\Psi_{\mu_n}|)|\Psi_{\mu_n}|)$, whose exponential order is
at least $\mu_1+\mu_n\ge\mu_{n+1}$ and whose polynomial weight is finite.
Increasing $A_{n+1}$ if necessary therefore establishes
\eqref{eq:expansion-induction}(b) at the next step, while condition~(a)
follows from the claim. This completes the induction.

\smallskip\noindent
\emph{Step 3. Angular degrees and conclusion.}
We now verify the angular-degree bound by induction on $j$. A homogeneous
mode at $\mu_j=\gamma_k$ has degree
$k=\mathfrak d(\gamma_k)\le M^*(\mu_j)$, while the pressure forcing at
$\mu_j=\rho$, when present, is radial. Since $\mathcal L$ preserves
exponential order, $\mathcal LS_{j-1}$ has no term at $\mu_j$. Thus every
other contribution to $\mathscr S_{\mu_j}$ is a product of $m\ge2$
previously constructed terms. If their exponents are
$\nu_1,\ldots,\nu_m$, with $\nu_1+\cdots+\nu_m=\mu_j$, the induction
hypothesis bounds their degrees by $M^*(\nu_1),\ldots,M^*(\nu_m)$.
For each factor choose a generator representation attaining the
corresponding maximum in $M^*$. Adding these representations gives
$\mathbf n\in\mathbb N_0^{(\mathcal G_F)}$ with
$|\mathbf n|\ge m\ge2$ and $\langle\mathbf n\rangle=\mu_j$.
Lemma~\ref{lem:harmonics} therefore bounds the degree of the product by
\[
 \sum_{\ell=1}^m M^*(\nu_\ell)
 =\sum_{\gamma\in\mathcal G_F}n_\gamma\,\mathfrak d(\gamma)
 \le\widetilde M(\mu_j)\le M^*(\mu_j).
\]
Solving the modal equations does not change the spherical degree, so
every spherical harmonic in $\Psi_{\mu_j}$ has degree at most $M^*(\mu_j)$.

At each fixed order, only finitely many terms occur. Products add finitely
many powers of $t$, a resonant modal solve raises the degree by at most
one, and Lemma~\ref{lem:linear-extraction} gives a finite remainder weight.
For each $j$, choose $M_j\in\mathbb N_0$ such that
\[
 M_j\ge A_j,\qquad
 M_j\ge\deg_t\!\bigl(e^{-\mu_jt}\Psi_{\mu_j}\bigr),
\]
where the degree of the zero polynomial is taken to be zero. Padding
missing powers with zero coefficients, write
\[
 \Psi_{\mu_j}(t,\theta)
 =e^{\mu_jt}\sum_{i=0}^{M_j}t^i C_{ji}(\theta).
\]
Each $C_{ji}$ is a finite sum of spherical harmonics of degree at most
$M^*(\mu_j)$. Hence $z=S_n+\mathcal R_{n+1}$ gives, for every $n\ge1$,
\[
 z(t,\theta)=\sum_{j=1}^{n}\sum_{i=0}^{M_j}
 C_{ji}(\theta)t^ie^{\mu_jt}
 +O\!\left((1+|t|)^{M_{n+1}}e^{\mu_{n+1}t}\right).
\]
Subsequent steps leave the previously constructed $\Psi_{\mu_j}$
unchanged, so these coefficients are the same for all truncation orders.
This proves the proposition.
\end{proof}

\begin{remark}
When $F\ne0$, the particular solution $\eta_\star$ used in the proof of
Proposition~\ref{prop:decay} contributes the known degree-zero term at
exponent $\rho$. If $\rho=\gamma_j$ for some $j\ge k_*$, this radial
term and the free degree-$j$ homogeneous term have the same exponent;
any nonlinear contribution at that order is determined by preceding
induction steps.
\end{remark}
\section{Existence and completion of the proofs}
Throughout this section we assume \eqref{H} and retain the notation
of Section~3.
The linear theory follows \cite{13z} and \cite[Section~IV]{8z}; we record only
the modifications caused by the H\'enon parameter and the modal cutoff.

For $i\in\mathbb N_0$, $\beta\in(0,1)$ and $\mu\in\mathbb R$, set
\[
\begin{aligned}
\|v\|_{C_\mu^i}
&=\sum_{j=0}^i\sup_{(-\infty,t_0]\times S^{N-1}}
e^{-\mu t}|\nabla^jv(t,\theta)|,\\
\|v\|_{C_\mu^{i,\beta}}
&=\|v\|_{C_\mu^i}
+\sup_{t\le t_0-1}e^{-\mu t}
[\nabla^iv]_{C^\beta([t-1,t+1]\times S^{N-1})},
\end{aligned}
\]
where the underlying cylinder is omitted from the notation when no confusion
can arise.  The corresponding weighted H\"older space consists of the
functions with finite norm.

Recall that $\gamma_k=\sigma_1^{(k)}>0$ for $k\ge k_*$.  Define
\[
\mathscr D_+=\{\gamma_k:k\ge k_*\}\subset\mathcal G_F\subset\mathcal I .
\]
If $\mu>0$ and $\mu\notin\mathscr D_+$, put
\begin{equation}\label{eq:weighted-cutoff}
\mathcal K(\mu)=\max\bigl(\{0,\ldots,k_*-1\}
\cup\{k\ge k_*:\gamma_k<\mu\}\bigr)
\end{equation}
(the second set may be empty), and let
\[
\mathcal H_{\le \mathcal K}=\bigoplus_{k=0}^{\mathcal K}\operatorname{span}
\{Q_1^k,\ldots,Q_{m_k}^k\}.
\]
Thus the full multiplicity of each spherical eigenspace lies on one side of
the weight.  For
$g\in C_\mu^{0,\beta}((-\infty,t_0]\times S^{N-1})$ we consider
\begin{equation}\label{eq:eqe}
\mathcal Lv=g.
\end{equation}

\begin{lem}\label{lem:schauder}
Let $\beta\in(0,1)$ and $\mu>0$.  Suppose that
$g\in C_\mu^{0,\beta}$, $\psi\in C^{2,\beta}(S^{N-1})$, and $v$ solves
\begin{equation}\label{eq:initial-problem}
\begin{cases}
\mathcal Lv=g&\text{in }(-\infty,t_0)\times S^{N-1},\\
v=\psi&\text{on }\{t_0\}\times S^{N-1},
\end{cases}
\end{equation}
with finite $C_\mu^0$ norm and local $C^{2,\beta}$ regularity.  Then
$v\in C_\mu^{2,\beta}$ and
\[
\|v\|_{C_\mu^{2,\beta}}
\le C\bigl(\|v\|_{C_\mu^0}+\|g\|_{C_\mu^{0,\beta}}
+e^{-\mu t_0}\|\psi\|_{C^{2,\beta}(S^{N-1})}\bigr),
\]
where $C$ is independent of $t_0$.
\end{lem}

\begin{proof}
For $t\le t_0-2$, apply the interior Schauder estimate on
$[t-2,t+2]\times S^{N-1}$; for $t_0-2\le t\le t_0$, apply the boundary
estimate on $[t_0-4,t_0]\times S^{N-1}$.  On every such cylinder,
$e^{-\mu s}$ and $e^{-\mu t}$ are uniformly comparable.  Multiplication by
$e^{-\mu t}$ and taking the supremum yields the asserted estimate.  This is
the standard weighted argument; compare \cite[Lemma~2.5]{13z} and
\cite[Lemma~13]{8z}.
\end{proof}

\begin{lem}\label{lem:apriori}
Let $\mu>0$, $\mu\notin\mathscr D_+$ and $\mathcal K=\mathcal K(\mu)$.  Suppose
$T\le t_0-4$, $g\in C^0([T,t_0]\times S^{N-1})$, and
$v\in C^2([T,t_0]\times S^{N-1})$ satisfies
\[
\begin{cases}
\mathcal Lv=g&\text{in }(T,t_0)\times S^{N-1},\\
v=0&\text{on }(\{T\}\cup\{t_0\})\times S^{N-1}.
\end{cases}
\]
If, for all $t\in[T,t_0]$, $0\le k\le \mathcal K$ and $1\le j\le m_k$,
\[
\int_{S^{N-1}}v(t,\theta)Q_j^k(\theta)\,\mathrm{d}\theta=0,
\]
then
\[
\sup_{[T,t_0]\times S^{N-1}}e^{-\mu t}|v|
\le C\sup_{[T,t_0]\times S^{N-1}}e^{-\mu t}|g|,
\]
where $C$ is independent of $T$ and $t_0$.
\end{lem}

\begin{proof}
The proof is the $t\to-\infty$ counterpart of
\cite[Lemma~2.6]{13z}; see also \cite[Lemma~14]{8z}.  We include the modal
point needed for the present cutoff.  If the estimate failed, there would be
sequences $T_i\le t_{0,i}-4$ and corresponding solutions $v_i$ with
right-hand sides $g_i$ such that
\[
 \|v_i\|_{C_\mu^0([T_i,t_{0,i}]\times S^{N-1})}=1,
 \qquad
 \|g_i\|_{C_\mu^0([T_i,t_{0,i}]\times S^{N-1})}\longrightarrow0.
\]
Choose $t_i^*\in(T_i,t_{0,i})$ at which the weighted maximum is attained and
set
\[
 \widetilde v_i(t,\theta)=e^{-\mu t_i^*}v_i(t+t_i^*,\theta),
 \qquad
 \widetilde g_i(t,\theta)=e^{-\mu t_i^*}g_i(t+t_i^*,\theta).
\]
Then $\widetilde v_i$ is defined on
$(T_i-t_i^*,t_{0,i}-t_i^*)$ and satisfies
\[
\sup_{S^{N-1}}|\widetilde v_i(0,\cdot)|=1,\qquad
|e^{-\mu t}\widetilde v_i(t,\theta)|\le1,\qquad
\mathcal L\widetilde v_i=\widetilde g_i\longrightarrow0
\]
locally uniformly.  Boundary $W^{2,q}$ estimates show that the translated
endpoints stay a fixed distance from $0$; at every finite limiting endpoint
the limit has zero boundary value.  Interior compactness therefore gives a
nonzero solution $\widehat v$ of $\mathcal L\widehat v=0$ on an interval
$(\tau_-,\tau_+)$ containing $0$, with
$|\widehat v|\le e^{\mu t}$.

The orthogonality assumptions kill every modal coefficient of $\widehat v$
with $k\le \mathcal K$.  If $k>\mathcal K$, then $k\ge k_*$ and
\[
\widehat v_j^k(t)=\kappa_1e^{\gamma_kt}+\kappa_2e^{\sigma_2^{(k)}t},
\qquad \gamma_k>\mu,\quad \sigma_2^{(k)}<0.
\]
At an infinite endpoint the weighted bound eliminates the corresponding
growing term, while at a finite endpoint the zero boundary value supplies the
remaining condition.  If both endpoints are finite, the determinant of the
two boundary conditions is nonzero because
$\gamma_k-\sigma_2^{(k)}>0$.  Hence $\kappa_1=\kappa_2=0$ in every case.
Completeness of the spherical harmonics gives $\widehat v=0$, a
contradiction.
\end{proof}

\begin{lem}\label{lem:low-modes}
Let $\beta\in(0,1)$, $\mu>0$, $\mu\notin\mathscr D_+$ and $\mathcal K=\mathcal K(\mu)$.
If $g\in C_\mu^{0,\beta}$ and $g(t,\cdot)\in \mathcal H_{\le \mathcal K}$ for every
$t\le t_0$, then \eqref{eq:eqe} has a unique solution
$v\in C_\mu^{2,\beta}$ satisfying $v(t,\cdot)\in \mathcal H_{\le \mathcal K}$.  The map
$g\mapsto v$ is linear and
\[
\|v\|_{C_\mu^{2,\beta}}\le C\|g\|_{C_\mu^{0,\beta}},
\]
with $C$ independent of $t_0$.
\end{lem}

\begin{proof}
Project onto $Q_j^k$, $0\le k\le \mathcal K$, and write the resulting equation as
$\mathcal L_kv_j^k=g_j^k$.  For distinct real roots,
\[
v_j^k(t)=\frac1{\sigma_1^{(k)}-\sigma_2^{(k)}}
\int_{-\infty}^t
\bigl(e^{\sigma_1^{(k)}(t-s)}-e^{\sigma_2^{(k)}(t-s)}\bigr)
g_j^k(s)\,\mathrm{d}s .
\]
The conjugate-root and repeated-root kernels are obtained by the usual
limits.  For $k\le \mathcal K$, every characteristic exponent has real part strictly
smaller than $\mu$; hence these integrals converge and give
$\|v_j^k\|_{C_\mu^0}\le C\|g_j^k\|_{C_\mu^0}$.  The one-dimensional local
Schauder estimate, together with the comparability of the weights on unit
intervals, yields the $C_\mu^{2,\beta}$ bound.  Summation over the finitely
many pairs $(k,j)$ proves existence and the estimate.  The same root
inequality excludes every nonzero homogeneous solution in this weighted
finite-dimensional class, proving uniqueness.  Compare
\cite[Lemma~2.7]{13z} and \cite[Lemma~15]{8z}.
\end{proof}

\begin{lem}\label{lem:high-modes}
Let $\beta\in(0,1)$, $\mu>0$, $\mu\notin\mathscr D_+$ and $\mathcal K=\mathcal K(\mu)$.
Assume $g\in C_\mu^{0,\beta}$ is orthogonal to all $Q_j^k$ with
$0\le k\le \mathcal K$.  Then \eqref{eq:eqe} has a unique solution
$v\in C_\mu^{2,\beta}$ with $v=0$ on
$\{t_0\}\times S^{N-1}$, and
\[
\|v\|_{C_\mu^{2,\beta}}\le C\|g\|_{C_\mu^{0,\beta}},
\]
where $C$ is independent of $t_0$.
\end{lem}

\begin{proof}
Let $\mathcal X_{\mathcal K}$ be the orthogonal complement in $H^1(S^{N-1})$ of the modes
of degree at most $\mathcal K$.  Since $\gamma_{\mathcal K+1}>\mu>0$,
$P_{\mathcal K+1}(0)=B_0-\lambda_{\mathcal K+1}<0$, and hence
\[
\int_{S^{N-1}}\bigl(|\nabla_\theta\xi|^2-B_0\xi^2\bigr)
\ge\left(1-\frac{B_0}{\lambda_{\mathcal K+1}}\right)
\int_{S^{N-1}}|\nabla_\theta\xi|^2,\qquad \xi\in\mathcal X_{\mathcal K}.
\]
Thus, on every finite cylinder $[T,t_0]\times S^{N-1}$, the weighted
Dirichlet functional
\[
\int_T^{t_0}\!\int_{S^{N-1}}e^{A_0t}
\bigl((\partial_t\xi)^2+|\nabla_\theta\xi|^2-B_0\xi^2+2g\xi\bigr)
\]
is coercive on $H_0^1$ functions taking values in $\mathcal X_{\mathcal K}$.  Its minimizer
$v_T$ is a weak solution of the projected, hence the full, Dirichlet problem.
Standard elliptic boundary regularity gives
$v_T\in C^{2,\beta}([T,t_0]\times S^{N-1})$, so Lemma~\ref{lem:apriori}
yields a $C_\mu^0$ bound independent of $T$ and $t_0$.  Local interior and
boundary Schauder estimates, followed by a diagonal argument as
$T\to-\infty$, produce a solution on the half-cylinder.  The $C_\mu^0$
estimate passes to the limit, and Lemma~\ref{lem:schauder} with zero boundary
data gives the stated $C_\mu^{2,\beta}$ estimate.  For uniqueness, project a
homogeneous solution onto spherical
harmonics: the weighted condition eliminates the negative root, and the
boundary value at $t_0$ eliminates $e^{\gamma_kt}$ for $k>\mathcal K$; the modes
$k\le \mathcal K$ are excluded as in Lemma~\ref{lem:low-modes}.  This is the same
finite-cylinder construction as in \cite[Lemma~2.8]{13z} and
\cite[Lemma~17]{8z}.
\end{proof}

\begin{thm}\label{thm:right-inverse}
Let $\beta\in(0,1)$, $\mu>0$, $\mu\notin\mathscr D_+$, and
$g\in C_\mu^{0,\beta}((-\infty,t_0]\times S^{N-1})$.  Then
\eqref{eq:eqe} has a solution $v\in C_\mu^{2,\beta}$ satisfying
\[
\|v\|_{C_\mu^{2,\beta}}\le C\|g\|_{C_\mu^{0,\beta}},
\]
where $C$ is independent of $t_0$.  The correspondence $g\mapsto v$ may be
chosen linear.
\end{thm}

\begin{proof}
Let $\Pi_{\le \mathcal K}$ be the projection onto $\mathcal H_{\le \mathcal K}$.  Apply Lemma
\ref{lem:low-modes} to $\Pi_{\le \mathcal K}g$ and Lemma~\ref{lem:high-modes} to
$(I-\Pi_{\le \mathcal K})g$, and add the two solutions.
\end{proof}

\begin{remark}
We denote this bounded linear right inverse by
\[
\mathcal L^{-1}:C_\mu^{0,\beta}\longrightarrow C_\mu^{2,\beta}.
\]
Its operator norm is independent of $t_0$; the boundary condition built into
the high-mode component is the one in Lemma~\ref{lem:high-modes}.
\end{remark}

Define
\[
\mathcal F(z)=\mathcal Lz-f(z)-Fe^{\rho t}.
\]
Thus $\mathcal F(z)=0$ is equation \eqref{3}.  For functions of $t$ alone we use the
corresponding one-dimensional weighted spaces.

\begin{thm}\label{thm:fixed-point}
Let $\beta\in(0,1)$.

\emph{(1) Radial case.}
Let $\mu>0$.  Suppose that, for some $t_1\le0$,
$\widehat z\in C^{2,\beta}((-\infty,t_1])$ satisfies
\[
|\widehat z(t)|+|\widehat z'(t)|\longrightarrow0
\quad\text{as }t\to-\infty,
\qquad
\mathcal F(\widehat z)\in C_\mu^{0,\beta}((-\infty,t_1]).
\]
Then there are $t_0\le t_1$ and a solution
$z\in C^{2,\beta}((-\infty,t_0])$ of $\mathcal F(z)=0$ such that
\[
|z(t)-\widehat z(t)|\le Ce^{\mu t}.
\]

\emph{(2) Nonradial case.}
Let $\mu>0$, $\mu\notin\mathscr D_+$.  Suppose that
$\widehat z\in C^{2,\beta}((-\infty,t_1]\times S^{N-1})$ satisfies
\[
|\widehat z(t,\theta)|+|\nabla\widehat z(t,\theta)|
\longrightarrow0
\quad\text{uniformly in }\theta\text{ as }t\to-\infty,
\]
and
\begin{equation}\label{64q}
\mathcal F(\widehat z)\in
C_\mu^{0,\beta}((-\infty,t_1]\times S^{N-1}).
\end{equation}
Then there are $t_0\le t_1$ and
$z\in C^{2,\beta}((-\infty,t_0]\times S^{N-1})$ satisfying $\mathcal F(z)=0$ and
\[
|z(t,\theta)-\widehat z(t,\theta)|\le Ce^{\mu t}.
\]
With the equation parameters, $\mu$ and $\beta$ fixed, $t_0$ and $C$
can be chosen uniformly for a family of approximate solutions on a
common half-cylinder, provided that their residual norms in
$C_\mu^{0,\beta}$ have a common bound and that
$|\widehat z|+|\nabla\widehat z|\to0$ uniformly over the family and in
$\theta$ as $t\to-\infty$.
\end{thm}

\begin{proof}
We prove (2); the radial proof is identical, using the Green kernel for
$\mathcal L_0$.  Since the roots of $P_0$ have negative real part, this
one-dimensional right inverse is bounded
$C_\mu^{0,\beta}\to C_\mu^{2,\beta}$ for every $\mu>0$.  The argument below
is the standard correction of an approximate solution; compare
\cite[Theorem~3.1]{13z} and \cite[Theorem~20]{8z}.

We use the elementary embedding
\begin{equation}\label{eq:weighted-C1-embedding}
\|h\|_{C_\mu^{0,\beta}}\le C_{\mu,\beta}\|h\|_{C_\mu^1},
\end{equation}
which follows from the mean-value inequality on unit cylinders and local
comparability of the weights.  Write
\[
\mathcal F(\widehat z+\phi)
=\mathcal F(\widehat z)+\mathcal L\phi-\mathcal N(\phi),
\]
where
\[
\mathcal N(\phi)=\lambda\bigl[(\widehat z+\Lambda+\phi)^{-p}
-(\widehat z+\Lambda)^{-p}
+p\Lambda^{-(p+1)}\phi\bigr],
\]
and set
\[
\mathcal T(\phi)=
\mathcal L^{-1}\bigl(-\mathcal F(\widehat z)+\mathcal N(\phi)\bigr).
\]
Let $C_{\mathcal L}$ be the norm of the right inverse and choose
$\mathcal B>0$ with
$\mathcal B\ge2C_{\mathcal L}\|\mathcal F(\widehat z)\|_{C_\mu^{0,\beta}}$.
There is a nondecreasing $\omega(t)\to0$ as $t\to-\infty$ such that
$|\widehat z|+|\nabla\widehat z|\le\omega(t)$.  After decreasing $t_0$,
\[
\omega(t_0)+\mathcal Be^{\mu t_0}\le\frac{\Lambda}{2},
\]
so all negative powers below are uniformly defined on the ball
$\{\|\phi\|_{C_\mu^{2,\beta}}\le\mathcal B\}$.

Taylor's formula and \eqref{eq:weighted-C1-embedding} give, with
$\widetilde\omega(t_0)=\omega(t_0)+\mathcal Be^{\mu t_0}$,
\[
\begin{aligned}
\|\mathcal N(\phi)\|_{C_\mu^{0,\beta}}
&\le C\widetilde\omega(t_0)\mathcal B,\\
\|\mathcal N(\phi_1)-\mathcal N(\phi_2)\|_{C_\mu^{0,\beta}}
&\le C\widetilde\omega(t_0)
\|\phi_1-\phi_2\|_{C_\mu^{2,\beta}}.
\end{aligned}
\]
These estimates follow, for example, by writing
\[
\mathcal N(\phi)=
-\lambda p\,\phi\int_0^1
\bigl[(\widehat z+\Lambda+s\phi)^{-(p+1)}
-\Lambda^{-(p+1)}\bigr]\,\mathrm{d}s .
\]
Since the norm of $\mathcal L^{-1}$ is independent of $t_0$, we may
further decrease $t_0$ so that
$CC_{\mathcal L}\widetilde\omega(t_0)\le1/2$.
The preceding estimates and the choice of $\mathcal B$ then make
$\mathcal T$ a self-map and a contraction on the above ball.
Its fixed point $\phi$ satisfies
$\mathcal F(\widehat z+\phi)=0$ and $|\phi|\le\mathcal Be^{\mu t}$.

For a family satisfying the additional hypotheses, choose $\mathcal B$
and $\omega$ uniformly over the family, using the residual norms on the
common initial half-cylinder.  The same smallness conditions then hold
for one sufficiently negative $t_0$, giving a common cylinder and a
uniform error bound.
\end{proof}

The approximate solutions below are finite sums of
polynomial-exponential modes.  The following elementary fact handles the
resonant terms.

\begin{lem}[Resonant modes]\label{lem:resonance}
Let $\ell\ge0$, $\widetilde\rho>0$, and let $q$ be a polynomial of degree
at most $\mathcal J$.  If $P_\ell(\widetilde\rho)\ne0$, there is a polynomial $\mathcal Q$ of
degree at most $\mathcal J$ such that
\[
\mathcal L_\ell(\mathcal Q(t)e^{\widetilde\rho t})
=q(t)e^{\widetilde\rho t}.
\]
If $P_\ell(\widetilde\rho)=0$, the same conclusion holds with
$\deg \mathcal Q\le \mathcal J+1$; in this case $\widetilde\rho$ is a simple root.
\end{lem}

\begin{proof}
This follows from the polynomial coefficient calculation in the
one-step extraction argument in the proof of
Proposition~\ref{prop:expansion}, with $\nu=\widetilde\rho>0$.
\end{proof}

\begin{thm}\label{thm:approx}
Fix $\beta\in(0,1)$, and let
\[
\rho=\frac{2p-\alpha}{p+1},\qquad
d_F=\begin{cases}F/P_0(\rho),&F\ne0,\\0,&F=0,\end{cases}
\qquad
\eta_\star(t)=d_Fe^{\rho t}.
\]
If $F\ne0$, then $\rho>0$, $P_0(\rho)>0$ and
$\mathcal L_0\eta_\star=Fe^{\rho t}$; when $F=0$,
$\eta_\star\equiv0$.

\emph{(1) Radial case.}
Assume $F\ne0$ and $\mu>\rho$.  Choose $n_\mu\ge1$ such that
$(n_\mu+1)\rho>\mu$, let $d_1=d_F,d_2,\ldots,d_{n_\mu}$ be the
coefficients in the recurrence \eqref{eq:radial-recurrence}, and set
\[
\widehat z(t)=\sum_{\ell=1}^{n_\mu}d_\ell e^{\ell\rho t}.
\]
Then, on $(-\infty,t_0]$ for some $t_0<0$,
\[
\mathcal F(\widehat z)=O(e^{(n_\mu+1)\rho t}),
\]
and $\widehat z$ satisfies the hypotheses of
Theorem~\ref{thm:fixed-point}(1).

\emph{(2) Nonradial case.}
Choose
\[
\mu>\gamma_{k_*},\qquad
\mu>\rho\ \text{if }F\ne0,\qquad
\mu\notin\mathcal I .
\]
Such weights exist because the exponent set $\mathcal I$ is locally
finite; in particular, $\mu\notin\mathscr D_+$.
Let
\[
\varphi(t,\theta)=
\sum_{\substack{k\ge k_*:\ \gamma_k<\mu}}
\sum_{j=1}^{m_k}c_j^kQ_j^k(\theta)e^{\gamma_kt}
\]
be a finite sum with at least one nonzero coefficient, and let $k_0$ be
the smallest degree with a nonzero coefficient.  Set
$\rho_{\min}=\gamma_{k_0}$ if $F=0$, and
$\rho_{\min}=\min\{\rho,\gamma_{k_0}\}$ if $F\ne0$.
There exist $t_0<0$, nonnegative integers $\mathcal J(\widetilde\rho)$,
and a finite sum
\begin{equation}\label{eq:phitilde-structure}
\widetilde\varphi(t,\theta)=
\sum_{\substack{\widetilde\rho\in\mathcal I_{\widetilde\rho}\\
\widetilde\rho<\mu}}
\sum_{j=0}^{\mathcal J(\widetilde\rho)}
\Psi_{\widetilde\rho,j}(\theta)t^je^{\widetilde\rho t},
\end{equation}
where each $\Psi_{\widetilde\rho,j}$ is a finite sum of spherical harmonics
of degree at most $\widetilde M(\widetilde\rho)$, such that
\[
\widehat z=\eta_\star+\varphi+\widetilde\varphi
\]
satisfies
\[
|\widehat z|+|\nabla\widehat z|\longrightarrow0
\quad\text{uniformly in }\theta\text{ as }t\to-\infty,
\qquad
\mathcal F(\widehat z)\in
C_\mu^{0,\beta}((-\infty,t_0]\times S^{N-1}).
\]
Hence $\widehat z$ satisfies the hypotheses of
Theorem~\ref{thm:fixed-point}(2).  Moreover:
\begin{enumerate}[(i)]
\item If
$\mathcal I_\rho\cap\mathcal I_{\widetilde\rho}=\varnothing$, then
$\mathcal J(\widetilde\rho)=0$ for every $\widetilde\rho$.
\item Every nonzero term of $\widetilde\varphi$ has exponent at least
$2\rho_{\min}$, and every positive-degree component has exponent strictly
larger than $\gamma_{k_0}$; in particular, no term
$t^ae^{\gamma_{k_0}t}Q_b^{k_0}$ with $a\in\mathbb N_0$ and
$1\le b\le m_{k_0}$ occurs in $\widetilde\varphi$.
Thus the leading nonradial mode prescribed by $\varphi$ is unchanged by
the correction.
\item The coefficients of $\widetilde\varphi$ are polynomials in $d_F$
and the $c_j^k$.  With the equation parameters, $\mu$, $\beta$ and the
spherical-harmonic basis fixed, $t_0$ and a bound for
$\|\mathcal F(\widehat z)\|_{C_\mu^{0,\beta}}$ can be chosen independently
of the coefficients satisfying $\max_{k,j}|c_j^k|\le1$.
Moreover, $\widehat z$ and its first derivatives tend to zero uniformly
over this family and in $\theta$ as $t\to-\infty$.
\end{enumerate}
\end{thm}

\begin{proof}
We only prove (2).  The proof of (1) is similar and easier.

\smallskip
\emph{Nonradial case.}
Starting from $\eta_\star+\varphi$, with
$\mathcal L(\eta_\star+\varphi)=Fe^{\rho t}$, we construct
$\widetilde\varphi$ by successively cancelling all nonlinear residual terms
with exponent below $\mu$, so
$\widehat z=\eta_\star+\varphi+\widetilde\varphi$ satisfies the hypotheses of
Theorem~\ref{thm:fixed-point}(2). Following
\cite[Proposition~3.3]{13z} and \cite[Proposition~22]{8z}, we track
the pressure generator, angular degrees, and possible resonance.

\smallskip\noindent
\emph{Step 1. Preparation and Taylor truncation.}
Write
\[
\eta_\star+\varphi
=\sum_{\substack{\gamma\in\mathcal G_F\\\gamma<\mu}}
\Theta_\gamma(\theta)e^{\gamma t},
\]
where $\Theta_\gamma$ has spherical-harmonic degree at most
$\mathfrak d(\gamma)$ and includes the constant coefficient $d_F$ when
$F\ne0$ and $\gamma=\rho$.  Products of at least two such terms have
exponents in $\mathcal I_{\widetilde\rho}$.
Let
\[
\widetilde\rho_1<\cdots<\widetilde\rho_M<\mu
\]
be all elements of $\mathcal I_{\widetilde\rho}$ below $\mu$; $M=0$ is
allowed.  Choose an integer $J\ge2$ such that
$(J+1)\min\mathcal G_F>\mu$, hence also
$(J+1)\rho_{\min}>\mu$.  This choice of $J$ is independent of the
coefficients $c_j^k$.
For a finite partial sum $\Upsilon$, define the truncated residual
\[
\mathscr E_J(\Upsilon)
:=\mathcal L\Upsilon-\sum_{m=2}^{J}a_m(\eta_\star+\Upsilon)^m.
\]
Then
\[
\mathcal F(\eta_\star+\Upsilon)
=\mathscr E_J(\Upsilon)+\mathscr R_J(\eta_\star+\Upsilon),
\]
where $\mathscr R_J(y)=\sum_{m=2}^{J}a_my^m-f(y)$ is the Taylor remainder.
If every nonzero term of $\Upsilon$ has exponent at least $\rho_{\min}$,
Taylor's formula gives, on a sufficiently negative half-cylinder,
\[
|\mathscr R_J(\eta_\star+\Upsilon)|
+|\nabla\mathscr R_J(\eta_\star+\Upsilon)|
\le C(1+|t|)^A e^{(J+1)\rho_{\min}t}
\]
for some nonnegative integer $A$.  Thus only $\mathscr E_J(\Upsilon)$
needs to be considered at exponents below $\mu$.

\smallskip\noindent
\emph{Step 2. Recursive cancellation.}
Starting from $\Upsilon_0=\varphi$, we add corrections only at the
exponents $\widetilde\rho_1,\ldots,\widetilde\rho_M$.
For every such partial sum, $\mathscr E_J(\Upsilon)$ is a finite sum of
polynomial--exponential modes with exponents in
$\mathcal I_{\widetilde\rho}$: indeed, $\mathcal L\varphi=0$,
$\mathcal L$ preserves the exponent of each correction, and every
nonlinear term is a product of at least two factors.  Hence the above
list contains all exponents below $\mu$ that need to be cancelled.

\begin{samepage}
For $s=1,\ldots,M$, suppose that $\Upsilon_{s-1}$ has been constructed
so that:
\begin{enumerate}[(a)]
\item $\mathscr E_J(\Upsilon_{s-1})$ has no terms with exponents
$\widetilde\rho_1,\ldots,\widetilde\rho_{s-1}$;
\item every correction added at an exponent $\nu$ has angular degree at
most $\widetilde M(\nu)$;
\item every nonzero correction has exponent at least $2\rho_{\min}$.
\end{enumerate}
These conditions hold for $s=1$, when no correction has yet been added.
\end{samepage}

At exponent $\widetilde\rho_s$, the linear part
$\mathcal L\Upsilon_{s-1}$ contributes nothing, since
$\mathcal L\varphi=0$ and all previous corrections have smaller
exponents.

We next verify the angular-degree bound for the coefficient at
$\widetilde\rho_s$.  Consider a nonlinear product of $m\ge2$ factors with
exponents $\nu_1,\ldots,\nu_m$ summing to $\widetilde\rho_s$.
For each factor choose $\mathbf n^{(i)}\in\mathbb N_0^{(\mathcal G_F)}$
with $\nu_i=\langle\mathbf n^{(i)}\rangle$ so that its angular degree is
at most $\sum_\gamma n_\gamma^{(i)}\mathfrak d(\gamma)$.
For a factor from $\eta_\star+\varphi$, use its single generator, so
$|\mathbf n^{(i)}|=1$.  For a correction, use (b) and choose a maximizer
in \eqref{eq:Mtilde}, so $|\mathbf n^{(i)}|\ge2$.
Set $\mathbf n=\sum_{i=1}^{m}\mathbf n^{(i)}$.  Then
$|\mathbf n|\ge m\ge2$ and
$\langle\mathbf n\rangle=\widetilde\rho_s$.
By Lemma~\ref{lem:harmonics}, the product has angular degree at most
\[
\sum_{i=1}^{m}\sum_\gamma n_\gamma^{(i)}\mathfrak d(\gamma)
=\sum_\gamma n_\gamma\mathfrak d(\gamma)
\le\widetilde M(\widetilde\rho_s).
\]
Consequently, the term of $\mathscr E_J(\Upsilon_{s-1})$ at this exponent
has the form
\[
\sum_{\ell=0}^{\widetilde M(\widetilde\rho_s)}
\sum_{b=1}^{m_\ell}
q_{s\ell b}(t)Q_b^\ell(\theta)e^{\widetilde\rho_st},
\]
where the $q_{s\ell b}$ are polynomials.  Lemma~\ref{lem:resonance}
provides polynomials $\mathcal Q_{s\ell b}$ satisfying
\[
\mathcal L_\ell
\bigl(\mathcal Q_{s\ell b}(t)e^{\widetilde\rho_st}\bigr)
=-q_{s\ell b}(t)e^{\widetilde\rho_st}.
\]
In a resonant mode we choose $\mathcal Q_{s\ell b}(0)=0$.
This normalization excludes an additional homogeneous term, whose
coefficient is already represented in $\varphi$.  Define
\[
\Upsilon_s(t,\theta)=\Upsilon_{s-1}(t,\theta)
+\sum_{\ell=0}^{\widetilde M(\widetilde\rho_s)}
\sum_{b=1}^{m_\ell}
\mathcal Q_{s\ell b}(t)Q_b^\ell(\theta)e^{\widetilde\rho_st}.
\]
The linear contribution of this correction cancels the term at
$\widetilde\rho_s$.  Every new nonlinear term contains this correction
and at least one further positive-exponent factor, and hence has exponent
strictly larger than $\widetilde\rho_s$.  Thus (a) holds at the next step.
The modal solves preserve angular degree, giving (b).
A nonzero new correction arises from a product of at least two existing
factors, each with exponent at least $\rho_{\min}$, so (c) also holds.
The induction therefore continues through $s=M$.  Set
\[
\widetilde\varphi=\Upsilon_M-\varphi,
\qquad \widehat z=\eta_\star+\Upsilon_M.
\]
When $M=0$, this means $\widetilde\varphi=0$.
The construction and the angular-degree bound give exactly
\eqref{eq:phitilde-structure}.

\smallskip\noindent
\emph{Step 3. Residual estimate.}
All terms of $\mathscr E_J(\Upsilon_M)$ with exponent below $\mu$ have
now been cancelled.  Since $\mu\notin\mathcal I$ and
$\mathcal I$ is locally finite, the gap
\[
\delta_\mu=\min\{\nu-\mu:\nu\in\mathcal I,\ \nu>\mu\}
\]
is positive.  Every remaining term of $\mathscr E_J(\Upsilon_M)$ has exponent at least
$\mu+\delta_\mu$.  The same is true of the exponent
$(J+1)\rho_{\min}$ in the Taylor-remainder estimate, since it belongs to
$\mathcal I$ and is larger than $\mu$.  Enlarging $A$ if necessary, we
therefore obtain, uniformly in $\theta$,
\[
|\mathcal F(\widehat z)(t,\theta)|
+|\nabla\mathcal F(\widehat z)(t,\theta)|
\le C(1+|t|)^A e^{(\mu+\delta_\mu)t},\qquad t\le t_0<0.
\]
The positive gap absorbs the polynomial factor, so
$\mathcal F(\widehat z)\in C_\mu^1$ and
\eqref{eq:weighted-C1-embedding} yields
$\mathcal F(\widehat z)\in C_\mu^{0,\beta}$.
Also, $\widehat z$ is a finite sum of smooth modes with positive
exponents, and hence $|\widehat z|+|\nabla\widehat z|\to0$ uniformly in
$\theta$ as $t\to-\infty$.  This verifies all hypotheses of
Theorem~\ref{thm:fixed-point}(2).

\smallskip\noindent
\emph{Step 4. Additional properties.}
For (i), if
$\mathcal I_\rho\cap\mathcal I_{\widetilde\rho}=\varnothing$, then
$P_\ell(\widetilde\rho_s)\ne0$ for every $\ell,s$, so the polynomials
$q_{s\ell b}$ and $\mathcal Q_{s\ell b}$ may be taken constant.

For (ii), the lower bound $2\rho_{\min}$ was proved in Step~2.
A further induction along the same recursion shows that every
positive-degree correction has exponent strictly larger than
$\gamma_{k_0}$.  Indeed, a product has a positive-degree component only
if at least one factor is nonradial.  Initially, every such factor
comes from $\varphi$ and has exponent at least $\gamma_{k_0}$.
At each later step, a newly generated positive-degree term contains
either an original nonradial factor or a previously constructed
positive-degree correction, together with at least one further factor
of positive exponent.  Its exponent is therefore strictly larger than
$\gamma_{k_0}$, so no later correction can alter the leading $k_0$-mode.

For (iii), the number of recursive steps and the Taylor order $J$ are
fixed, and each resonant solve increases the polynomial degree by at
most one.  Hence all polynomial degrees have bounds independent of the
coefficients $c_j^k$.  On these fixed polynomial spaces, the modal
solves in Step~2 are fixed linear maps, and all other operations are
additions and multiplications.  Thus the coefficients of
$\widetilde\varphi$ are polynomials in $d_F$ and the $c_j^k$, and are
uniformly bounded on the compact set $\max_{k,j}|c_j^k|\le1$.
With the equation parameters, $\mu$, $\beta$ and the basis fixed, the
finite polynomial--exponential form therefore gives
\[
|\widehat z(t,\theta)|+|\nabla\widehat z(t,\theta)|
\le C(1+|t|)^A e^{\mu_1t},\qquad t\le t_0<0,
\]
where $\mu_1=\min\mathcal G_F>0$, and $C$ and $A$ are independent of the
coefficients.  This proves the uniform decay and allows one to choose
one $t_0$ with $|\widehat z|\le\Lambda/2$ for the whole family.
The finite exponent set and the gap $\delta_\mu$ are also independent
of the coefficients.  Using $\mu_1$ as a common lower exponent in the
Taylor-remainder estimate of Step~1, and then applying Step~3, gives a
uniform $C_\mu^{0,\beta}$ bound for the residual.  This proves (iii).
\end{proof}

\medskip
\noindent\textbf{Completion of the proofs of
Theorems~\ref{thm:radial} and~\ref{thm:nonradial}.}
With $d$ as in \eqref{qcc}, the inverse cylindrical transformation is
\begin{equation}\label{eq:reconstruct}
u(x)=\Lambda|x|^d+|x|^dz(\ln|x|,\theta),
\qquad \theta=\frac{x}{|x|}.
\end{equation}
If $\mathcal F(z)=0$ on $(-\infty,t_0]\times S^{N-1}$, then
\eqref{eq:reconstruct} solves \eqref{QQ} on
$B_{e^{t_0}}\setminus\{0\}$.  If also $z\to0$ uniformly in $\theta$
as $t\to-\infty$, decreasing $t_0$ gives $|z|\le\Lambda/2$.
Thus $u>0$ and, since $d>0$, $u(x)\to0$ as $x\to0$.

\begin{proof}[Proof of Theorem~\ref{thm:radial}]
Fix $\mu>\rho$.  Theorems~\ref{thm:approx}(1) and
\ref{thm:fixed-point}(1) yield a radial solution $z$ satisfying
\[
z(t)=\widehat z(t)+O(e^{\mu t})=O(e^{\rho t}).
\]
By \eqref{eq:reconstruct}, the corresponding function $u$ is a positive
radial rupture solution with
\[
 u(r)=\Lambda r^d\bigl(1+O(r^\rho)\bigr)
 \qquad\text{as }r\to0^+.
\]
Proposition~\ref{prop:radial-exp} then gives the asserted expansion for every
integer $K\ge1$ and $d_i=F^iD_i(\lambda,p,\alpha,N)$.
\end{proof}

\begin{proof}[Proof of Theorem~\ref{thm:nonradial}]
\emph{(i)} Fix a real-valued, $L^2(S^{N-1})$-normalized spherical
harmonic $Q=Q_1^{k_*}$ of degree $k_*$, and a weight
\[
 \mu>\gamma_{k_*},\qquad
 \mu>\rho\ \text{if }F\ne0,\qquad
 \mu\notin\mathcal I.
\]
For each $\mathfrak c\in(0,1]$, apply Theorem~\ref{thm:approx}(2) with
\[
 \varphi_{\mathfrak c}(t,\theta)=\mathfrak cQ(\theta)e^{\gamma_{k_*}t}.
\]
This gives an approximate solution
$\widehat z_{\mathfrak c}=\eta_\star+\varphi_{\mathfrak c}
+\widetilde\varphi_{\mathfrak c}$.
By Theorem~\ref{thm:approx}(2)(iii), these approximate solutions have a
common residual bound and tend to zero with their first derivatives
uniformly in $\mathfrak c$.  The uniform assertion in
Theorem~\ref{thm:fixed-point}(2) therefore gives solutions
\[
 z_{\mathfrak c}=\widehat z_{\mathfrak c}+O(e^{\mu t})
\]
on one common half-cylinder $(-\infty,t_0]\times S^{N-1}$, with the
error bound uniform for $0<\mathfrak c\le1$.
Decreasing this common $t_0$ once more gives
$|z_{\mathfrak c}|\le\Lambda/2$ for all $\mathfrak c$.
Thus \eqref{eq:reconstruct} yields positive rupture solutions
$u_{\mathfrak c}$ on the same punctured ball $B_R\setminus\{0\}$,
where $R=e^{t_0}$.  For this choice of $\varphi_{\mathfrak c}$,
$k_0=k_*$ and $\rho_{\min}=\mu_1$.
By Theorem~\ref{thm:approx}(2)(ii) and $\mu>\mu_1$,
\[
 u_{\mathfrak c}(x)=\Lambda|x|^d\bigl(1+O(|x|^{\mu_1})\bigr)
 \qquad\text{as }|x|\to0^+.
\]

It remains to distinguish these solutions.  The radial term $\eta_\star$
is orthogonal to $Q$.  By Theorem~\ref{thm:approx}(2)(ii), every
positive-degree component of $\widetilde\varphi_{\mathfrak c}$ has
exponent strictly larger than $\gamma_{k_*}$, while the error is
$O(e^{\mu t})$ with $\mu>\gamma_{k_*}$.  Since $Q$ is real-valued and
normalized,
\[
 e^{-\gamma_{k_*}t}
 \int_{S^{N-1}}z_{\mathfrak c}(t,\theta)Q(\theta)\,\mathrm{d}\theta
 \longrightarrow \mathfrak c
 \qquad\text{as }t\to-\infty.
\]
Thus every $u_{\mathfrak c}$ is nonradial, and
$u_{\mathfrak c_1}\ne u_{\mathfrak c_2}$ whenever
$\mathfrak c_1\ne\mathfrak c_2$.
The family $\{u_{\mathfrak c}:0<\mathfrak c\le1\}$ therefore consists
of infinitely many distinct positive nonradial asymptotically radial
rupture solutions on $B_R\setminus\{0\}$.

\smallskip
\emph{(ii)} Let $u$ be any positive asymptotically radial rupture
solution of \eqref{QQ} in the sense of \eqref{defn}, and set
\[
z(t,\theta)=r^{-d}u(r\theta)-\Lambda,\qquad t=\ln r.
\]
Proposition~\ref{prop:decay} gives $z=O(e^{\mu_1t})$, while
Proposition~\ref{prop:expansion} gives its expansion to every prescribed
finite order.  Substituting $t=\ln r$ in \eqref{eq:reconstruct} yields
the decay estimate and asymptotic expansion asserted in
Theorem~\ref{thm:nonradial}(ii).
\end{proof}

\section{Far-field classification and global radial continuation}\label{sec:infinity}

We prove a radial initial-value continuation lemma, then
Theorem~\ref{thm:infinity-classification}, and apply the lemma to local rupture
solutions in Corollary~\ref{cor:radial-global-connection}.
Assume \eqref{H} throughout; exterior results impose no condition at the origin.

\subsection{Radial continuation from positive initial data}

We consider the initial value problem
\eqref{eq:radial-initial-value-problem}.  Set
$a_\alpha(r):=\lambda r^\alpha$.  For use both here and in the
negative-pressure far-field argument, define, for $s>0$,
\begin{equation}\label{eq:infinity-potential}
 V(s):=\int_1^s(1-\xi^{-p})\,\mathrm{d}\xi=
 \begin{cases}
 \displaystyle s-1-\frac{s^{1-p}-1}{1-p},&p\ne1,\\[7pt]
 s-1-\log s,&p=1.
 \end{cases}
\end{equation}
Then $V'(s)=1-s^{-p}$, $V''(s)=ps^{-p-1}>0$, and $V\ge0$, with equality
only at $s=1$.  Moreover,
\begin{equation}\label{eq:infinity-potential-coercive}
 V(s)\longrightarrow\infty\quad(s\to\infty),\qquad
 \lim_{s\downarrow0}V(s)=\infty\quad\text{if }p\ge1.
\end{equation}

\begin{lem}[Radial continuation]\label{lem:radial-continuation}
For every $R_0>0$, $u_0>0$ and $u_1\in\mathbb R$, the initial value
problem \eqref{eq:radial-initial-value-problem} has a unique positive
solution on a maximal interval $[R_0,R^*)$, with
$R^*\in(R_0,\infty]$.
\begin{enumerate}[(i)]
\item If $R^*<\infty$, then
\[
 \sup_{R_0\le r<R^*}\bigl(u(r)+|u'(r)|\bigr)<\infty,
 \qquad \lim_{r\uparrow R^*}u(r)=0.
\]
\item If $p\ge1$, or if $F\ge0$ and $u_1\ge0$, then $R^*=\infty$.
\end{enumerate}
\end{lem}

\begin{proof}
The radial equation is equivalent to
\[
 \frac{\mathrm{d}}{\mathrm{d}r}
 \begin{pmatrix}u\\[1mm]u'\end{pmatrix}
 =\begin{pmatrix}u'\\[1mm]
 -\dfrac{N-1}{r}u'+a_\alpha(r)u^{-p}+F\end{pmatrix}.
\]
Its right-hand side is smooth for $r>0$ and $u>0$.  Local Picard existence
and uniqueness therefore give the unique maximal positive solution.

Suppose that $R^*<\infty$.  Since $a_\alpha$ is positive and continuous
on $[R_0,R^*]$, it is bounded above and away from zero there.  Define
\begin{equation}\label{eq:radial-continuation-energy}
 \mathcal E(r):=\frac12u'(r)^2+a_\alpha(r)V(u(r))\ge0.
\end{equation}
Using $a_\alpha'(r)=\alpha a_\alpha(r)/r$ and the radial equation gives
\begin{equation}\label{eq:radial-continuation-energy-derivative}
\begin{aligned}
 \mathcal E'(r)
 &=u'(r)\bigl[u''(r)+a_\alpha(r)(1-u(r)^{-p})\bigr]
   +a_\alpha'(r)V(u(r))\\
 &=-\frac{N-1}{r}u'(r)^2+(a_\alpha(r)+F)u'(r)
   +\frac{\alpha}{r}a_\alpha(r)V(u(r)).
\end{aligned}
\end{equation}
By Young's inequality,
\[
\begin{aligned}
 \mathcal E'(r)
 &\le\frac12(a_\alpha(r)+F)^2+\frac12u'(r)^2
     +\frac{|\alpha|}{R_0}a_\alpha(r)V(u(r))\\
 &\le\mathcal C\bigl(1+\mathcal E(r)\bigr),
\end{aligned}
\]
where $\mathcal C>0$ is independent of $r\in[R_0,R^*)$.
Gronwall's inequality yields
\[
 1+\mathcal E(r)\le
 \bigl(1+\mathcal E(R_0)\bigr)e^{\mathcal C(r-R_0)}
 \le\bigl(1+\mathcal E(R_0)\bigr)e^{\mathcal C(R^*-R_0)}.
\]
Thus $u'$ and $V(u)$ are uniformly bounded.  The first property in
\eqref{eq:infinity-potential-coercive} also bounds $u$ from above.
In particular, $u$ is Lipschitz on $[R_0,R^*)$ and has a finite limit
$u_*\ge0$ as $r\uparrow R^*$.

If $u_*>0$, then $u$ is bounded away from zero near $R^*$; the radial
equation bounds $u''$, so $u'$ is Lipschitz near $R^*$ with finite limit
$u_*'$. The first-order system's integral form permits adjoining the
endpoint data $(u(R^*),u'(R^*))=(u_*,u_*')$. Local Picard existence and
uniqueness then extend the solution beyond $R^*$, preserving positivity
by continuity and contradicting maximality. Thus $u_*=0$, proving (i).

If $p\ge1$, the limit $u(r)\to0$ at a finite endpoint is impossible,
because $V(u)$ is bounded and $V(s)\to\infty$ as $s\downarrow0$.
Hence $R^*=\infty$ in this case.  If instead $F\ge0$ and $u_1\ge0$, then
\[
 \bigl(r^{N-1}u'(r)\bigr)'
 =r^{N-1}\bigl(a_\alpha(r)u(r)^{-p}+F\bigr)>0
\]
throughout the maximal interval.  Thus $u'(r)>0$ for $r>R_0$ and
$u(r)\ge u_0>0$.  A finite maximal endpoint is again excluded by (i).
This proves (ii).
\end{proof}

\begin{remark}
Lemma~\ref{lem:radial-continuation} uses no condition at the origin; for
$0<p<1$ it does not cover arbitrary negative initial slopes because
$V(0^+)=p/(1-p)<\infty$.
\end{remark}

\subsection{Far-field classification}

We first record the elementary one-dimensional fact used in the
non-radial argument.  For $y\in C^2((R_0,\infty))$, write
\[
 L_Ny:=\bigl(r^{N-1}y'(r)\bigr)'.
\]

\begin{lem}\label{lem:infinity-integration}
Let $y\in C^2((R_0,\infty))$; all $o(1)$ symbols below refer to
$r\to\infty$.
\begin{enumerate}[(i)]
\item If
\[
 L_Ny=r^{N-1}\bigl(C_*+o(1)\bigr),\qquad C_*\ne0,
\]
then
\[
 y(r)=\frac{C_*}{2N}r^2+o(r^2).
\]

\item If $C_*>0$, $\zeta>-2$, and
\[
 L_Ny=C_*r^{N-1+\zeta}\bigl(1+o(1)\bigr),
\]
then
\begin{equation}\label{eq:infinity-integration-positive}
 y(r)=\frac{C_*}{(N+\zeta)(\zeta+2)}
 r^{\zeta+2}\bigl(1+o(1)\bigr).
\end{equation}

\item If $C_*>0$, $\zeta\le-2$, and
\[
 L_Ny=C_*r^{N-1+\zeta}\bigl(1+o(1)\bigr),
\]
then $y$ is not asymptotic to a positive
multiple of a positive power of $r$.
\end{enumerate}
\end{lem}

\begin{proof}
Fix $R_1>R_0$.  For $\nu>-1$,
\[
 \int_{R_1}^r s^\nu\bigl(1+o(1)\bigr)\,\mathrm{d}s
 =\frac{r^{\nu+1}}{\nu+1}\bigl(1+o(1)\bigr).
\]
Every integration below starts at $R_1$.
Applying this identity proves (i) and (ii) easily. 

For (iii), first let $\zeta=-2$.  If $N\ge3$, two integrations give
\[
 y(r)=\frac{C_*}{N-2}\log r+o(\log r)+C_0+C_1r^{2-N},
\]
whereas for $N=2$ they give
\[
 y(r)=\frac{C_*}{2}(\log r)^2+o\bigl((\log r)^2\bigr)
      +C_1\log r+C_0.
\]
Neither expression is asymptotic to a positive power of $r$.

Now let $\zeta<-2$.  If $N=2$, the right-hand side defining $L_Ny$
is integrable at infinity; hence $ry'(r)$ has a finite limit and
$y(r)=O(1+\log r)$.  If $N\ge3$, then $y'(r)=O(r^{\zeta+1})$ when
$-N<\zeta<-2$, $y'(r)=O(r^{1-N}\log r)$ when $\zeta=-N$, and
$r^{N-1}y'(r)$ has a finite limit when $\zeta<-N$.  In all three cases
$y$ remains bounded.  This proves (iii).
\end{proof}

\begin{proof}[Proof of Theorem~\ref{thm:infinity-classification}]
\emph{(1) Radial case.}
We establish the asymptotic formulas for an arbitrary positive radial solution $u$ of
\eqref{eq:infinity-equation} on $\Omega_{R_0}$.

\medskip\noindent
\emph{The case $F>0$.}
Set
\[
 h(r):=\lambda r^\alpha u(r)^{-p}\ge0.
\]
The radial equation takes the flux form
\begin{equation}\label{eq:infinity-radial-flux}
 \bigl(r^{N-1}u'(r)\bigr)'=r^{N-1}\bigl(F+h(r)\bigr).
\end{equation}
Fix $R_1>R_0$.  Integration from $R_1$ to $r$ gives
\begin{equation}\label{eq:infinity-Fpositive-uprime}
 u'(r)=\frac{F}{N}r+c_0r^{1-N}
 +r^{1-N}\int_{R_1}^r s^{N-1}h(s)\,\mathrm{d}s,
\end{equation}
where
\[
 c_0=R_1^{N-1}u'(R_1)-\frac{F}{N}R_1^N.
\]
The last term in \eqref{eq:infinity-Fpositive-uprime} is nonnegative and
$c_0r^{1-N}=o(r)$.  Hence, for all sufficiently large $r$,
\[
 u'(r)\ge\frac{F}{2N}r.
\]
This yields
$u(r)\ge cr^2$ for some $c>0$.  Since $\alpha<2p$,
\[
 0\le h(r)\le\lambda c^{-p}r^{\alpha-2p}=o(1).
\]
It follows that
\[
 r^{-N}\int_{R_1}^r s^{N-1}h(s)\,\mathrm{d}s\longrightarrow0.
\]
Equation \eqref{eq:infinity-Fpositive-uprime} therefore gives
$u'(r)=\frac{F}{N}r+o(r)$.  One more integration proves
the first line of \eqref{eq:infinity-radial-conclusion}.

\medskip\noindent
\emph{The case $F=0$.}
Introduce the Emden--Fowler variables
\[
 t=\log r,\qquad u(r)=\Lambda r^{d}U(t),
\]
and recall $A_0=N-2+2d>0$ from \eqref{eq:AB}; set $c_d:=d(d+N-2)>0$, so that
$\Lambda^{p+1}=\lambda/c_d$ by \eqref{qcc}.  Using $d(p+1)=\alpha+2$, the
radial equation \eqref{eq:infinity-equation} with $F=0$ is equivalent to the
autonomous equation
\begin{equation}\label{eq:exterior-autonomous}
 U_{tt}+A_0U_t+c_d\bigl(U-U^{-p}\bigr)=0,\qquad U>0,
\end{equation}
and the balanced profile corresponds to the equilibrium $U\equiv1$.  Define,
for $U>0$,
\begin{equation}\label{eq:exterior-potential}
\Phi(U):=c_d\int_1^U\bigl(s-s^{-p}\bigr)\,\mathrm{d}s=
c_d\begin{cases}
 \dfrac{U^2-1}{2}-\dfrac{U^{1-p}-1}{1-p}, & p\ne1,\\[8pt]
 \dfrac{U^2-1}{2}-\log U, & p=1.
\end{cases}
\end{equation}
Since $s-s^{-p}$ is negative on $(0,1)$ and positive on $(1,\infty)$, we have
$\Phi\ge0$, with $\Phi(U)=0$ if and only if $U=1$.  For a solution of
\eqref{eq:exterior-autonomous}, set
$E(t)=\tfrac12U_t(t)^2+\Phi(U(t))$; then
\begin{equation}\label{eq:exterior-energy-identity}
 E'(t)=U_t\bigl[U_{tt}+c_d\bigl(U-U^{-p}\bigr)\bigr]=-A_0U_t^2\le0.
\end{equation}
Let $E_\infty:=\lim_{t\to\infty}E(t)\ge0$.  Since $u$ is only defined on the
open set $\Omega_{R_0}$, fix $R_1>R_0$ and set
$t_{R_1}:=\log R_1$; the estimates below concern $t\ge t_{R_1}$.

To prove $U(t)\to1$, it suffices to show that the orbit $(U,U_t)$
eventually remains in a compact subset of $(0,\infty)\times\mathbb R$.
Indeed, by \cite[Lemmas~6.5--6.6]{Teschl2012}, its $\omega$-limit set
is nonempty, compact, and invariant. By
\cite[Theorem~6.15]{Teschl2012}, the energy
$\tfrac12 U_t^2+\Phi(U)$ is constant on this set, so
\eqref{eq:exterior-energy-identity} implies that $U_t=0$ there.
Invariance and \eqref{eq:exterior-autonomous} then show that this
limit set consists only of the equilibrium $(1,0)$.
Consequently, \cite[Lemma~6.7]{Teschl2012} yields
$(U(t),U_t(t))\to(1,0)$ and $E_\infty=0$.

Suppose first that $p\ge1$.  Then $\Phi(U)\to\infty$ both as $U\to0^+$
and as $U\to\infty$, so the bound $E(t)\le E(t_{R_1})$ confines $(U,U_t)$
to a compact subset of $(0,\infty)\times\mathbb R$.  The preceding
argument gives $U(t)\to1$, which proves
the second line of \eqref{eq:infinity-radial-conclusion} in this case.

Suppose now that $0<p<1$, so that the barrier
$\Phi_0:=\Phi(0^+)=c_d(1+p)/[2(1-p)]$ is finite, while $\Phi(U)\to\infty$
as $U\to\infty$ still holds.  Thus the energy bound controls $U$ from
above and $|U_t|$, but by itself does not keep $U$ away from zero.  We
will prove $E_\infty<\Phi_0$, which supplies the missing positive lower
bound and allows us to apply the same compactness argument.

Indeed, if $E_\infty<\Phi_0$, then for all sufficiently large $t$,
\[
 \Phi(U(t))\le E(t)\le\frac{E_\infty+\Phi_0}{2}<\Phi_0.
\]
Since $\Phi(U)\to\Phi_0$ as $U\downarrow0$, this strict sublevel bound
keeps $U$ uniformly away from zero. Together with the bounds on $U$
and $U_t$, it gives the compactness required above, hence $U(t)\to1$.

To establish the strict inequality, suppose instead that
$E_\infty\ge\Phi_0$.  Then $E(t)\ge\Phi_0$ for every $t\ge t_{R_1}$.

If $U(t^*)\in(0,1)$ for some $t^*\ge t_{R_1}$, then, as long as
$U\in(0,1)$,
\[
 \tfrac12U_t^2=E-\Phi(U)\ge\Phi_0-\Phi(U)>0,
\]
so $U_t$ cannot vanish and keeps a fixed sign there.  A zero of $U_t$
at $U=1$ is also impossible, since it would give
$E=\Phi(1)=0<\Phi_0$.

If this sign is negative, then $U$ decreases strictly.  For $0<U<1$,
\[
 \Phi_0-\Phi(U)
 =c_d\left(\frac{U^{1-p}}{1-p}-\frac{U^2}{2}\right)
 \ge\Phi_0 U^{1-p},
\]
where we used $U^2\le U^{1-p}$.  Hence the time of descent to zero
is at most
\[
 \int_0^{U(t^*)}\frac{\mathrm{d}U}{\sqrt{2\bigl(\Phi_0-\Phi(U)\bigr)}}
 \le\frac1{\sqrt{2\Phi_0}}
 \int_0^{U(t^*)}U^{-(1-p)/2}\,\mathrm{d}U<\infty.
\]
The last integral is finite because $(1-p)/2<1$.  Thus $U$ approaches
zero at a finite time, contradicting positivity of the given exterior
solution for all $r>R_0$.

If the sign is positive, set $U_*:=U(t^*)\in(0,1)$.  Since $\Phi$ is
strictly decreasing on $(0,1)$, until $U$ first reaches $1$ we have
\[
 U_t(t)\ge\sqrt{2\bigl(\Phi_0-\Phi(U(t))\bigr)}
 \ge\sqrt{2\bigl(\Phi_0-\Phi(U_*)\bigr)}=: c_*>0.
\]
Integration gives $U(t)\ge U_*+c_*(t-t^*)$, so $U$ reaches $1$ in finite
time and crosses it with positive velocity.  Any later return to
$(0,1)$ must occur with $U_t<0$, since zero velocity at $U=1$ has been
excluded, and would lead to the finite-time contradiction above.

Under the assumption $E_\infty\ge\Phi_0$, global positivity therefore
forces $U(t)\ge1$ for all sufficiently large $t$.  Together with the
energy bound, this confines the tail of the orbit $(U,U_t)$ to a compact
subset of $(0,\infty)\times\mathbb R$.  LaSalle's invariance
principle then yields $E_\infty=0$, contradicting $E_\infty\ge\Phi_0>0$.
Hence $E_\infty<\Phi_0$, and the sublevel argument above gives
$U(t)\to1$, proving the second line of
\eqref{eq:infinity-radial-conclusion} for $0<p<1$ as well.

\medskip\noindent
\emph{The case $F<0$.}
Write $F=-G$, where $G>0$, and set
\begin{equation}\label{eq:infinity-Fnegative-parameters}
 m_F:=\frac{\alpha}{p},\qquad
 C_F:=\left(\frac{\lambda}{G}\right)^{1/p},\qquad
 \kappa:=\frac{2-m_F}{2}.
\end{equation}
Here $p\ge1$ and $0\le\alpha<2p$, so $0\le m_F<2$ and $\kappa>0$.  Introduce
\begin{equation}\label{eq:infinity-Fnegative-scaling}
 u(r)=C_Fr^{m_F} y(r),\qquad
 s=\frac{1}{\kappa}\sqrt{\frac{G}{C_F}}\,r^\kappa,
 \qquad Y(s)=y(r).
\end{equation}
Since $\mathrm{d}s/\mathrm{d}r=\kappa s/r$, direct differentiation gives
\[
\begin{aligned}
 u'(r)&=C_Fr^{m_F-1}\bigl(m_FY+\kappa sY'\bigr),\\
 u''(r)&=C_Fr^{m_F-2}\bigl[\kappa^2s^2Y''
 +\kappa(2m_F+\kappa-1)sY' +m_F(m_F-1)Y\bigr].
\end{aligned}
\]
Consequently,
\[
\begin{split}
 u''+\frac{N-1}{r}u'=C_Fr^{m_F-2}\bigl[&\kappa^2s^2Y''
 +\kappa(2m_F+\kappa+N-2)sY'\\
 &+m_F(m_F+N-2)Y\bigr].
\end{split}
\]
The choices of $m_F$ and $C_F$ imply
$\lambda r^\alpha u^{-p}=GY^{-p}$, while
$C_Fr^{m_F-2}=G/(\kappa^2s^2)$.  Thus $Y$ satisfies
\begin{equation}\label{eq:infinity-Y-equation}
 Y''+\frac{\Gamma_F}{s}Y'+\frac{h_0}{s^2}Y=Y^{-p}-1,
\end{equation}
where
\begin{equation}\label{eq:infinity-Gamma-H}
 \Gamma_F=\frac{2N-2+3m_F}{2-m_F}>0,\qquad
 h_0=\frac{4m_F(m_F+N-2)}{(2-m_F)^2}\ge0.
\end{equation}

Use the potential $V$ defined in \eqref{eq:infinity-potential}.
Since $p\ge1$, it tends to infinity both at zero and at infinity by
\eqref{eq:infinity-potential-coercive}.  Define
\begin{equation}\label{eq:infinity-energy}
 E_F(s):=\frac12Y'(s)^2+V(Y(s))+\frac{h_0}{2s^2}Y(s)^2.
\end{equation}
Using \eqref{eq:infinity-Y-equation}, we obtain
\begin{equation}\label{eq:infinity-energy-derivative}
 E_F'(s)=-\frac{\Gamma_F}{s}Y'(s)^2
 -\frac{h_0}{s^3}Y(s)^2\le0.
\end{equation}
Hence $E_F(s)\downarrow E_{F,\infty}$ for some $E_{F,\infty}\ge0$.

For sufficiently large $S$, the bound $E_F(s)\le E_F(S)$ and
\eqref{eq:infinity-potential-coercive} give constants $0<Y_-<Y_+$ and
$M>0$ such that
\begin{equation}\label{eq:infinity-Y-bounds}
 Y_-\le Y(s)\le Y_+,\qquad
 |Y'(s)|+|Y''(s)|\le M,\qquad s\ge S.
\end{equation}
We claim that $E_{F,\infty}=0$.  Suppose otherwise and put
\[
 D(s):=\int_s^{s+1}Y'(\xi)^2\,\mathrm{d}\xi.
\]
First, $\liminf_{s\to\infty}D(s)>0$.  Indeed, if $D(s_n)\to0$ along
some sequence $s_n\to\infty$, the uniform bound for $Y''$ in
\eqref{eq:infinity-Y-bounds} implies
\[
 \sup_{[s_n,s_n+1]}|Y'|\longrightarrow0.
\]
After passing to a subsequence, $Y(s_n)\to L\in[Y_-,Y_+]$, and hence
$Y_n(\tau):=Y(s_n+\tau)\to L$ in $C^1([0,1])$.  Integrating
\eqref{eq:infinity-Y-equation} over $[s_n,s_n+1]$ gives
\[
\begin{split}
 &Y'(s_n+1)-Y'(s_n)
 +\int_0^1\frac{\Gamma_F}{s_n+\tau}Y'(s_n+\tau)\,\mathrm{d}\tau\\
 &\qquad
 +\int_0^1\frac{h_0}{(s_n+\tau)^2}Y(s_n+\tau)\,\mathrm{d}\tau
 =\int_0^1\bigl(Y(s_n+\tau)^{-p}-1\bigr)\,\mathrm{d}\tau .
\end{split}
\]
The left-hand side tends to zero by
$\sup_{[s_n,s_n+1]}|Y'|\to0$ and
\eqref{eq:infinity-Y-bounds}, whereas the right-hand side tends to
$L^{-p}-1$.  Thus $L^{-p}-1=0$, so $L=1$.  Evaluating
\eqref{eq:infinity-energy} at $s=s_n$ now yields $E_F(s_n)\to0$,
contradicting $E_{F,\infty}>0$.

Thus $D(s)\ge d_0>0$ for all sufficiently large $s$.  Integration of
\eqref{eq:infinity-energy-derivative} over $[s,s+1]$ gives
\[
 E_F(s+1)-E_F(s)\le-\frac{\Gamma_F d_0}{s+1}.
\]
Summing this inequality along consecutive unit intervals contradicts
$E_F\ge0$, because the harmonic series diverges.  Therefore $E_{F,\infty}=0$.
Every term in \eqref{eq:infinity-energy} is nonnegative, so
$V(Y(s))\to0$ and hence $Y(s)\to1$.  Reversing
\eqref{eq:infinity-Fnegative-scaling} proves
the third line of \eqref{eq:infinity-radial-conclusion}.

\medskip\noindent
\emph{(2) Possibly non-radial case.}
Assume
\eqref{eq:infinity-power-hypothesis}.  Define the normalized spherical
averages
\[
 \overline u(r):=\frac{1}{|S^{N-1}|}\int_{S^{N-1}}u(r,\theta)\,\mathrm{d}\theta,
 \qquad
 \overline{u^{-p}}(r):=\frac{1}{|S^{N-1}|}\int_{S^{N-1}}u(r,\theta)^{-p}\,\mathrm{d}\theta.
\]
The uniformity in \eqref{eq:infinity-power-hypothesis} gives
\begin{equation}\label{eq:infinity-averages-asymptotic}
 \overline u(r)=A_\infty r^{\sigma_\infty}\bigl(1+o(1)\bigr),
 \qquad
 \overline{u^{-p}}(r)=A_\infty^{-p}r^{-p\sigma_\infty}\bigl(1+o(1)\bigr).
\end{equation}

The spherical average of the angular Laplacian vanishes.  Therefore
averaging \eqref{eq:infinity-equation} gives the exact identity
\begin{equation}\label{eq:infinity-averaged-equation}
\begin{split}
 \bigl(r^{N-1}\overline u'(r)\bigr)'
 &=r^{N-1}\bigl[\lambda r^\alpha \overline{u^{-p}}(r)+F\bigr]\\
 &=r^{N-1}\bigl[a_\infty r^{\eta_\infty}(1+o(1))+F\bigr],
\end{split}
\end{equation}
where
\[
 a_\infty:=\lambda A_\infty^{-p}>0,
 \qquad \eta_\infty:=\alpha-p\sigma_\infty.
\]

Suppose first that $F>0$.  If $\eta_\infty<0$, the bracket in
\eqref{eq:infinity-averaged-equation} is $F+o(1)$, and
Lemma~\ref{lem:infinity-integration}(i) gives
\[
 \overline u(r)=\frac{F}{2N}r^2+o(r^2).
\]
Comparison with \eqref{eq:infinity-averages-asymptotic} yields
$\sigma_\infty=2$ and $A_\infty=F/(2N)$.  If $\eta_\infty=0$, the same lemma, now with
coefficient $a_\infty+F>0$, again forces $\sigma_\infty=2$; but then $\eta_\infty=0$ would
give $\alpha=2p$, contrary to \eqref{H}.  If $\eta_\infty>0$, the nonlinear
term dominates $F$, and Lemma~\ref{lem:infinity-integration}(ii) gives
$\sigma_\infty=\eta_\infty+2$.  Hence
\[
 \sigma_\infty=\alpha-p\sigma_\infty+2
 \quad\Longrightarrow\quad
 \sigma_\infty=d=\frac{\alpha+2}{p+1}.
\]
Since $\eta_\infty=\sigma_\infty-2>0$, this would imply $d>2$, equivalently
$\alpha>2p$, again a contradiction.  Thus the first line of
\eqref{eq:infinity-nonradial-conclusion} follows.

Now suppose that $F=-G<0$, where $G>0$.  If $\eta_\infty<0$, the bracket in
\eqref{eq:infinity-averaged-equation} is $-G+o(1)$, so
Lemma~\ref{lem:infinity-integration}(i) gives
\[
 \overline u(r)=-\frac{G}{2N}r^2+o(r^2),
\]
contradicting $\overline u(r)>0$.  If $\eta_\infty>0$, the nonlinear term dominates
$-G$ and the argument above yields $\sigma_\infty=d>2$, contradicting
$\alpha<2p$.  Hence $\eta_\infty=0$, and
\[
 \bigl(r^{N-1}\overline u'(r)\bigr)'
 =r^{N-1}\bigl[a_\infty-G+o(1)\bigr].
\]
If $a_\infty-G<0$, Lemma~\ref{lem:infinity-integration}(i) contradicts the
positivity of $\overline u$.  If $a_\infty-G>0$, the same lemma forces $\sigma_\infty=2$,
whereas $\eta_\infty=0$ gives $\sigma_\infty=\alpha/p<2$.  Therefore $a_\infty=G$, and
\[
 \sigma_\infty=\frac{\alpha}{p},
 \qquad
 A_\infty=\left(\frac{\lambda}{-F}\right)^{1/p}.
\]
This proves the third line of
\eqref{eq:infinity-nonradial-conclusion}.

Finally, let $F=0$.  If $\sigma_\infty\le0$, then
$\eta_\infty=\alpha-p\sigma_\infty\ge\alpha>-2$.  Lemma~\ref{lem:infinity-integration}(ii)
makes $\overline u(r)$ asymptotic to a positive multiple of $r^{\eta_\infty+2}$, whose
exponent is positive, contradicting
\eqref{eq:infinity-averages-asymptotic}.  Hence $\sigma_\infty>0$.  If
$\eta_\infty\le-2$, Lemma~\ref{lem:infinity-integration}(iii) contradicts this
positive-power asymptotic.  Thus $\eta_\infty>-2$, and
Lemma~\ref{lem:infinity-integration}(ii), together with
\eqref{eq:infinity-averages-asymptotic}, gives
\[
\begin{aligned}
 \sigma_\infty&=\eta_\infty+2=\alpha-p\sigma_\infty+2,\\
 A_\infty&=\frac{\lambda A_\infty^{-p}}{(N+\eta_\infty)(\eta_\infty+2)}
  =\frac{\lambda A_\infty^{-p}}{(\sigma_\infty+N-2)\sigma_\infty}.
\end{aligned}
\]
Consequently,
\[
 \sigma_\infty=d=\frac{\alpha+2}{p+1},
 \qquad
 A_\infty^{p+1}=\frac{\lambda}{d(d+N-2)},
\]
so $A_\infty=\Lambda$.  This proves the second line of
\eqref{eq:infinity-nonradial-conclusion}.  The remainders remain uniform
in $\theta$.
\end{proof}

\begin{remark}
In the radial $F<0$ argument, $\alpha<2p$, $\alpha\ge0$, and $p\ge1$
are used respectively for $\kappa>0$, $h_0\ge0$, and coercivity of $V$ at zero;
$\alpha\ge0$ is not needed for continuation.
\end{remark}

\subsection{Global continuation of radial rupture solutions}

\begin{proof}[Proof of Corollary~\ref{cor:radial-global-connection}]
Fix $r_0\in(0,R)$.  If $F>0$, then $u(r)\to0$ as $r\downarrow0$ and
$u(r_0)>0$.  The mean value theorem gives a point $r_1\in(0,r_0)$ with
$u'(r_1)>0$.  Since
\[
 \bigl(r^{N-1}u'(r)\bigr)'
 =r^{N-1}\bigl(\lambda r^\alpha u(r)^{-p}+F\bigr)>0,
\]
we have $u'(r_0)>0$.  Lemma~\ref{lem:radial-continuation}(ii), applied
with initial radius $r_0$, therefore extends $u$ to a positive solution
on $[r_0,\infty)$ for every $p>0$.

If $F<0$ and $p\ge1$, the same lemma applies to the data
$(u(r_0),u'(r_0))$ without any restriction on the sign of $u'(r_0)$ or
of $\alpha$.  It again gives a positive continuation on $[r_0,\infty)$.
In both cases, local uniqueness makes the continuation agree with the
given solution on their common interval.  Combining it with the original
solution on $(0,r_0)$ gives the unique positive radial continuation to
$\mathbb R^N\setminus\{0\}$.

The far-field assertions follow from the $F>0$ and $F<0$ cases of
Theorem~\ref{thm:infinity-classification}(1), respectively;
the latter requires the additional condition $\alpha\ge0$.
Finally, Theorem~\ref{thm:radial} supplies the local radial rupture
solutions to which this argument applies.
\end{proof}

\begin{acknowledgments}
Yanyan Zhang is supported by the National Natural Science Foundation of China
(Grant No. 12271505) and the Science and Technology Commission of Shanghai
Municipality (Grant No. 22DZ2229014).
\end{acknowledgments}

\section*{Author Declarations}

\subsection*{Conflict of Interest}
The authors have no conflicts to disclose.

\subsection*{Author Contributions}
\noindent\textbf{Yunxiao Li:} Writing -- original draft (equal).
\textbf{Yanyan Zhang:} Writing -- original draft (equal).

\section*{Data Availability}
Data sharing is not applicable to this article as no new data were created or analyzed in this study.

\renewcommand{\refname}{References}
\bibliography{references}

\end{document}